\documentclass[11pt]{article}
\usepackage{amssymb}

\newcommand{\Z}{\mathbb{Z}}
\newcommand{\F}{\mathbb{F}}
\newcommand{\PP}{{\rm Prob}}

\newcommand{\Em}[1]{\textbf{#1}}

\newtheorem{theorem}{Theorem}[section]
\newtheorem{conjecture}[theorem]{Conjecture}
\newtheorem{proposition}[theorem]{Proposition}
\newtheorem{lemma}[theorem]{Lemma}
\newtheorem{corollary}[theorem]{Corollary}
\newtheorem{defn}[theorem]{Definition}

\newtheorem{unremarks}[theorem]{Remarks}

\def\ee{\epsilon}

\def\pp{{\bf p}}

\def\P{{\mathbb P}}
\def\E{{\mathbb E}}
\def\Var{{\rm Var}\,}
\def\Cov{{\rm Cov}\,}

\def\Cox{\hfill \Box}

\def\Z{{\mathbb Z}}

\def\R{{\mathbb R}}

\def\F{{\cal F}}

\def\disp{\displaystyle}

\def\hyp{{\cal H}}
\def\supp{{\rm supp}\,}
\def\cemetery{\Delta}
\def\graph{{\cal G}}

\def\sector{W}
\def\symdif{{\oplus}}
\def\tree{{\bf T}}
\def\marked{U}
\def\vmarked{\tilde{U}}
\def\cT{{\cal W}}
\def\VS{{\bf V}}
\def\FF{{\mathbb F}}
\def\A{{\cal A}}
\def\B{{\cal B}}
\def\SS{{\cal S}}

\begin{document}

\setcounter{page}{0}

\centerline{\huge \bf Sharp Transitions in Making Squares}

\vskip1in

\centerline{\Large Ernie Croot\footnote{Supported in part by an NSF
grant.}}

\centerline{\large School of Mathematics} \centerline{\large
Georgia Tech} \centerline{\large Atlanta, GA 30332-0160}
\centerline{\large ecroot@math.gatech.edu}
\bigskip

\bigskip

\centerline{\Large Andrew Granville \footnote{Partiellement
soutenu par une bourse de la Conseil de recherches en sciences
naturelles et en g\'enie du Canada.}}

\centerline{\large D\'epartment de math\'ematiques et de
statistique} \centerline{\large Universit\'e de Montr\'eal}
\centerline{\large Montr\'eal QC H3C 3J7, Canada}
\centerline{\large andrew@dms.umontreal.ca}
\bigskip
\bigskip

\centerline{\Large Robin Pemantle\footnote{Supported in part by NSF Grant DMS-01-03635.}}

\centerline{\large Department of Mathematics}
\centerline{\large University of Pennsylvania}
\centerline{\large 209 S. 33rd Street}
\centerline{\large Philadelphia, Pennsylvania 19104, USA}
\centerline{\large pemantle@math.upenn.edu}

\bigskip

\bigskip

\centerline{\Large Prasad Tetali\footnote{Supported in part by NSF Grants DMS-0401239 and DMS-0701043.}}

\centerline{\large School of Mathematics and School of Computer Science}
\centerline{\large Georgia Tech} \centerline{\large Atlanta, GA
30332-0160} \centerline{\large tetali@math.gatech.edu}

\newpage

\begin{abstract}

In many integer factoring algorithms, one produces a sequence of
integers (created in a pseudo-random way), and wishes to rapidly
determine a subsequence whose product is a square (which we call a
{\sl square product}).  In his lecture at the 1994 International
Congress of Mathematicians, Pomerance observed that the following
problem encapsulates all of the key issues:\ {\it Select integers
$a_1,a_2,\dots,$ at random from the interval $[1,x]$, until some
(non-empty) subsequence has product equal to a square. Find  good
estimate for the expected stopping time of this process.}
A good solution to this problem should
help one to determine the optimal choice of parameters for one's
factoring algorithm, and therefore this is a central question.

\

Pomerance (1994), using an idea of Schroeppel (1985), showed that
with probability $1-o(1)$ the first subsequence whose product
equals a square occurs after at least $J_0^{1-o(1)}$ integers have
been selected, but no more than $J_0$, for an appropriate
(explicitly determined) $J_0=J_0(x)$. Herein we determine this
expected stopping time up to a constant factor, tightening
Pomerance's interval to
$$
[ (\pi/4)(e^{-\gamma} - o(1))J_0,\ (e^{-\gamma} + o(1)) J_0],
$$
where $\gamma = 0.577...$ is the Euler-Mascheroni constant.
We will also
confirm the well established belief that, typically, none of the
integers in the square product have large prime factors.

We
believe that there should, in fact, be a sharp threshold for this
stopping time, that it should occur with probability $1-o(1)$
after at least $\{e^{-\gamma}-o(1)\} J_0$ integers have been
selected, but no more than $\{e^{-\gamma}+o(1)\} J_0$, with $J_0(x)$
as before.

Our proofs use  methods somewhat different from previous articles on this
subject. The heart of the proof of the upper bound lies in delicate
calculations in probabilistic graph theory, supported by
comparative estimates on smooth numbers using precise information on saddle points.

\end{abstract}

\newpage
\section{Introduction}

Several algorithms for factoring integers $n$ (including Dixon's
random squares algorithm \cite{dixon}, the quadratic sieve
\cite{pom0}, the multiple polynomial quadratic sieve
\cite{silverman}, and the number field sieve \cite{buhler} -- see
\cite{pom3} for a nice expository article on factoring algorithms)
work by generating a pseudorandom sequence of integers
$a_1,a_2,...$, with each
$$
a_i\ \equiv\ b_i^2 \pmod{n},
$$
until some subsequence of the $a_i$'s has product equal to a
square. Say we have such a subsequence
$$
a_{i_1}, ..., a_{i_k},\ {\rm where}\ Y^2\ =\ a_{i_1}\cdots
a_{i_k},
$$
and set
$$
X^2\ =\ (b_{i_1}\cdots b_{i_k})^2.
$$
Then
$$
n\ |\ Y^2 - X^2\ =\ (Y-X)(Y+X),
$$
and there is a fair chance that ${\rm gcd}(n, Y - X)$ is a
non-trivial factor of $n$.  If so, we have factored $n$.

In his lecture at the 1994 International Congress of
Mathematicians, Pomerance \cite{pom1,pom2} observed that in the
(heuristic) analysis of such factoring algorithms one assumes that
the pseudo-random sequence $a_1,a_2,...$ is close enough to random
that we can make predictions based on this assumption. Hence it
makes sense to formulate this question in its own right, in
particular to determine whether this part of the factoring
algorithm can be significantly sped up.
\bigskip

\noindent {\bf Pomerance's Problem.} Select positive integers
$a_1,a_2,\dots \leq x$ independently at random (that is, $a_j=m$
with probability $1/x$ for each integer $m,\ 1\leq m\leq x$),
until some subsequence of the $a_i$'s has product equal to a
square. When this occurs, we say that the sequence has a {\it
square dependence}. What is the expected stopping time of this
process ? \medskip

To discuss the history of this problem, and our own work, we need
to introduce some notation:\ Let $\pi(y)$ denote the number of
primes up to $y$. Call $n$ a $y$-{\sl smooth integer} if all of
its prime factors are $\leq y$, and let $\Psi(x,y)$ denote the
number of $y$-smooth integers up to $x$.  Let $y_0=y_0(x)$ be a
value of $y$ which maximizes $\Psi(x,y)/{y}$, and let
\begin{eqnarray} \label{J0_def}
J_0(x)\ :=\   \frac{ \pi(y_0)}{\Psi(x,y_0)} \cdot x.
\end{eqnarray}
In Pomerance's problem, let $T$ be the smallest integer $t$ for
which $a_1,...,a_t$ has a square dependence (note that $T$ is
itself a random variable). In 1985, Schroeppel gave a simple
argument to justify that for any $\epsilon>0$ we have
$$
\PP( T\ <\ (1+\epsilon)J_0(x))\ =\ 1-o(1)
$$
as $x\to \infty$, and in 1994 Pomerance showed that
$$
\PP(T\ >\ J_0(x)^{1-\epsilon})\ =\ 1-o(1).
$$
as $x\to \infty$. Therefore there is a transition from ``unlikely
to have a square product'' to ``almost certain to have a square
product'' at $T=J_0(x)^{1+o(1)}$. Pomerance asked in [3] whether
there is a sharper transition, and we conjecture that $T$ has a
{\it sharp threshold}: This would mean that there exists a
function $f(x)$ such that for every $\epsilon > 0$,
\begin{equation} \label{Tf}
\PP(T \in [ (1-\epsilon)f(x),\ (1+\epsilon)f(x)])\ =\ 1 - o(1)
\end{equation}
as $x\to \infty$. In fact we believe that this threshold is
$f(x)=e^{-\gamma} J_0(x)$:
\bigskip
\begin{conjecture}\label{conjecture}  For every $\epsilon > 0$ we have
\begin{equation} \label{Tg}
\PP(T \in [ (e^{-\gamma}-\epsilon) J_0(x),\ (e^{-\gamma}+\epsilon)
J_0(x)])\ =\ 1 - o(1),
\end{equation}
as $x\to \infty$, where $\gamma = 0.577...$ is the
Euler-Mascheroni constant.
\end{conjecture}

The constant $e^{-\gamma}$ in this conjecture is well-known to
number theorists. It appears as the ratio of the proportion of
integers free of prime divisors smaller than $y$, to the
proportion of integers up to $y$ that are prime. However this is not
how it appears in our discussion, and we have failed to find a more
direct route to this prediction.
\bigskip

The bulk of this article will be devoted to establishing the upper bound
in the above conjecture.
We will prove something a little weaker than the conjectured lower bound:

\begin{theorem}\label{main_theorem}  We have
$$
\PP(T\ \in\ [ (\pi/4)(e^{-\gamma} - \epsilon) J_0(x),\ (e^{-\gamma} + \epsilon) J_0(x)])\
=\ 1 - o(1),
$$
for any $\epsilon > 0$ as $x\to\infty$.
\end{theorem}

To obtain the lower bound in our theorem, we obtain a good upper
bound on the expected number of sub-products of the large prime
factors of the $a_i$'s that equal a square,  which allows us to
bound the probability that such a sub-product exists, for
$T<(\pi/4)(e^{-\gamma} - o(1)) J_0(x)$. This is the ``first moment
method''.

Schroeppel established his upper bound, $T\leq (1+o(1))J_0(x)$, by
showing that by then one expects more than $\pi(y_0)$ $y_0$-smooth
integers amongst $a_1,a_2,\dots ,a_T$, which guarantees that the
sequence has a square dependence. (To see this, create a matrix
over $\mathbb F_2$ whose columns are indexed by the primes up to
$y_0$, whose rows are indexed by the numbers $i$ such that $a_i$ is
$y_0$-smooth, and whose $(i,p)$th entry is given by the exponent on $p$
in the factorization of $a_i$, for each $y_0$-smooth $a_i$. Then a
square dependence amongst the $a_i$ is equivalent to a dependence
amongst the corresponding rows of our matrix, so that we are
guaranteed a square dependence once the matrix has more than
$\pi(y_0)$ rows.)  If we replace the complicated random model which creates
this matrix by one in which any given row appears as a
row of this matrix with equal probability then one  expects  a
linear dependence only once the matrix has more than $\pi(y_0)-O(1)$ rows
(see section 3.1 of \cite{CGPT} for details; also see \cite{Calkin} for a lower bound
in a related model of choosing binary vectors of fixed weight  randomly, until finding
a $GF(2)$-dependent set).

Schroeppel's approach is not only good for theoretical analysis,
in practice one searches among the $a_i$ for $y_0$-smooth integers
and hunts amongst these for a square dependence, using linear
algebra in $\mathbb F_2$ on the primes' exponents. Computing
specialists have also found that it is easy and profitable to keep
track of $a_i$ of the form $s_i q_i$, where $s_i$ is $y_0$-smooth
and $q_i$ is a prime exceeding $y_0$; if both $a_i$ and $a_j$ have
exactly the same large prime factor $q_i=q_j$ then their product
is a $y_0$-smooth integer times a square, and so can be used in
our matrix as an extra smooth number. This is called the {\sl
large prime variation}, and the upper bound in Theorem 1 of \cite{CGPT} is
obtained by computing the limit of this method  (to obtain
a constant, in place of $e^{-\gamma}$ which is a tiny bit smaller than $3/4$).

One can also consider the {\sl double large prime variation} in
which one allows two largish prime factors so that, for example,
the product of three $a_i$s of the form $pqs_1, prs_2, qrs_3$ can
be used as an extra smooth number. Experience has shown that each
of these variations has allowed a small speed up of various
factoring algorithms (though at the cost of some non-trivial extra
programming), and a long open question has been to formulate all
of the possibilities for multi-large prime variations and to
analyze how they affect the running time. Sorting out this
combinatorial maze has been the most difficult part of our work. 

When our process terminates (at time $T$) we have some subset $I$
of $a_1,...,a_T$, including $a_T$, whose product equals a
square.\footnote{Note that $I$ is unique, else if we have two such
subsets $I$ and $J$ then $(I\cup J)\setminus (I\cap J)$ is also a
set whose product equals a square, but does not contain $a_T$, and
so the process would have stopped earlier than at time $T$.}
It is not hard to show that this square product is $T^2$-smooth 
(see Section~3.2 of~\cite{CGPT}); here we give a more precise 
idea of what $I$ looks like:

\begin{theorem} \label{Theorem_2} \

a)\  In the special case that for  $\epsilon > 0$, conditional on the event $\{T < (\pi/4)(e^{-\gamma} - \epsilon)
J_0(x)\}$, we find that $I$ consists of a single number $a_i$ (which
is therefore a square) with probability $1-o(1)$.
\medskip

b) \ In general, with probability $1-o(1)$, we have that
\begin{equation} \label{1.1}
y_0 \exp( -(c_3+\epsilon) \sqrt{\log y_0})\ \leq\ |I|\ \leq\ y_0
\exp( (c_3+\epsilon) \sqrt{\log y_0})],
\end{equation}
where $c_3=\sqrt{2-\log 2}$. In other words, when the algorithm
terminates the square product $I$ is, almost certainly, composed
of $y_0^{1+o(1)}=J_0(x)^{1/2+o(1)}$ numbers $a_i$.
\medskip

c) \ Also, with probability $1-o(1)$ all the elements of $I$ are
$$
y_0^2 \exp(  (2+\epsilon) \sqrt{\log y_0\log\log y_0}) {\rm -smooth}.
$$
\end{theorem}

The last part of this result confirms the long held suspicion that
the earliest occurring square products are almost always composed
only of smooth numbers with a suitable smoothness parameter,
though the smoothness bound that we give may be significantly
larger than is possible, for all we know.

We expect that one can give more precise descriptions of $I$,
specifying more precisely how large $I$ is, and improving the
smoothness bound on the elements of $I$, perhaps even to
$y_0\phi(x)$ for any function $\phi$ for which $\phi(x)\to\infty$
as $x\to\infty$.
\bigskip

There are now several theorems along the lines of Conjecture 1 in
the literature, including some quite general approaches.
Friedgut's theorem \cite{friedgut}, characterizing a {\em coarse
threshold} for monotone or symmetric\footnote{That is, invariant
under permutations of the elements involved.} graph properties,
has been instrumental in proving the existence of a sharp
threshold for several graph properties. However it does not seem
to be applicable in the present context, since the square
dependence problem is not symmetric. Bourgain's strengthening of
sorts of Friedgut's theorem (see the appendix to \cite{friedgut})
is in principle applicable in the present context, though various
researchers have not yet succeeded in doing so.
\bigskip

Pomerance's main goal in enunciating the random squares problem was to provide a model that would prove useful in analyzing the running time of factoring algorithms, such as the quadratic sieve.
In \cite{CGPT} we analyzed the running time of Pomerance's random squares problem
to show  that the running time
will be inevitably dominated by  finding the actual
square product once we have enough integers. Indeed this carries
over to an analysis of the  quadratic sieve factoring algorithm
(and presumably the other factoring algorithms as well); a consequence is that to
optimize the running time of the quadratic sieve we
look for a square dependence among the $y$-smooth integers with $y$  significantly smaller than $y_0$, so that   Pomerance's problem is not
quite so germane to the question as it had at first appeared. Anyway, see
\cite{CGPT} for further discussion of these issues.
\bigskip

The paper is organized as follows. In section 2, we derive the
necessary technical lemmas involving smooth numbers. In section 3,
we   derive the lower bound for $T$
given in
Theorem~\ref{main_theorem}, and develop these ideas to prove
Theorem~\ref{Theorem_2}. Finally, in section 4, we develop our analysis of
multiprime variations.

\section{Smooth numbers}

In previous analyses of these questions, authors have typically
used estimates for $\Psi(x,y)$ for $y$ a fixed power of $y_0$. In
this range one can determine an asymptotic for $\Psi(x,y)$ in
terms of a saddle point, an implicit quantity. It has proved to be
difficult to deduce an asymptotic for $\Psi(x,y)$, or even
something close, in terms of simple explicit functions. One of the
key innovations in this article is to by-pass this issue by
comparing values of $\Psi(x,y)$ for different, but closely
related, values of $x$ and $y$: Since the saddle points are not
too different one can obtain sharp explicit estimates for the
ratio of two such $\Psi$-values.  In this technical section we
deduce several such results, primarily from the deep work of
Hildebrand and Tenenbaum \cite{hild}, which will come in useful
later.

\subsection{Classical smooth number estimates}

 From \cite{hild} we have that the estimate
\begin{equation} \label{2.1}
\Psi(x,y) = x\rho(u) \left\{ 1 + O\left( \frac{\log (u+1)}{\log y} \right)
\right\} \quad {\rm\ as\ }\quad x \rightarrow \infty
\quad {\rm where}\quad x=y^u,
\end{equation}
holds in the range
\begin{equation} \label{2.2}
\exp \left ( (\log\log x)^2 \right )\ \leq\ y\ \leq\ x,
\end{equation}
where $\rho(u)=1$ for $0 \le u \le 1$, and where
$$
\rho(u) = \frac{1}{u} \int_{u-1}^{u} \rho(t)\, dt \ \ {\rm for\ all}\
u>1.
$$
This function $\rho(u)$ satisfies
$$
\rho(u)\ =\ \exp ( -(u + o(u)) \log u );
$$
and so
\begin{equation} \label{psi_estimate}
\Psi(x,y)\ =\ x \exp(-(u+o(u))\log u).
\end{equation}

Now let
$$
L:=L(x) = \exp \left( \sqrt{ \frac 12 \log x \log\log x} \right).
$$
Then, using (\ref{psi_estimate}) we deduce that for $\beta > 0$,
\begin{equation} \label{psi_estimate_2}
\Psi(x, L(x)^{\beta+o(1)})\ =\ x L(x)^{-1/\beta + o(1)}.
\end{equation}
 From this one can easily deduce that
\begin{equation} \label{J0y0}
y_0(x)=L(x)^{1+o(1)},\ {\rm and\ }
J_0(x)=y_0^{2-\{1+o(1)\}/\log\log y_0}=L(x)^{2+o(1)},
\end{equation}
where $y_0$ and $J_0$ are as in the introduction (see
(\ref{J0_def})). From this we can deduce the following basic
estimate, which we will use in later proofs:

\begin{lemma} \label{first_lemma}
Fix constant $\beta > 0$. If $y = y_0^{\beta + o(1)}$ then
$$\\
{\Psi(x,y)/y \over \Psi(x,y_0)/y_0}\ =\ y_0^{2 - \beta - \beta^{-1} + o(1)}.
$$
\end{lemma}

\subsection{Hildebrand-Tenenbaum saddle point method estimates}

\noindent
For any $\alpha>0$, one has
\begin{equation} \label{2.3}
\Psi(x,y)\leq \sum_{n\leq x\atop P(n)\leq y} (x/n)^\alpha \leq x^{\alpha}\xi(\alpha,y),
\end{equation}
where
$$
\xi(s,y)\ =\ \prod_{p\le y}\Bigl(1-\frac{1}{p^s}\Bigr)^{-1}.
$$
Define $\alpha=\alpha(x,y)$ to be the solution to
\begin{equation}\label{alpha_def}
 \log x = \sum_{p\leq y} \frac{\log p}{p^\alpha-1} .
\end{equation}
By \cite[Theorem 1 and (7.19)]{hild} we obtain in the range (\ref{2.2}) with $u\to \infty$,
\begin{equation} \label{2.4}
\Psi(x,y)\sim \frac{x^{\alpha}\xi(\alpha,y)}{\alpha\sqrt{2\pi\log x\log y}}.
\end{equation}

Let $\xi=\xi(u)$ be the solution to $e^\xi=u\xi+1$
so that
\begin{equation} \label{2.4b}
\xi(u)\ =\ \log(u\log u) + {(1+o(1))\log\log u \over \log u},\ {\rm as\ }u \to \infty.
\end{equation}
Note also that $\xi'(u)\sim 1/u$.
In the range (\ref{2.2}) it turns out that
\begin{equation} \label{2.5}
(1-\alpha(x,y)) \log y = \xi(u) +O(1/u)
\end{equation}
which implies that
\begin{equation} \label{2.6}
y^{1-\alpha} = e^{\xi(u)} (1+O(1/u)) = u\xi(u)(1+O(1/u)).
\end{equation}
So, for
$$
y\ =\ L(x)^{\beta+o(1)}\ =\ y_0^{\beta + o(1)}
$$
we have
\begin{equation} \label{y_alpha}
y^{1-\alpha}\ \sim\ \beta^{-2} \log y\ \sim\ \beta^{-1} \log y_0.
\end{equation}

\noindent
By \cite[Theorem 3]{hild} and (\ref{2.5}) above, we have
\begin{equation} \label{2.7}
\Psi\left( \frac xd,y\right) = \frac 1{d^{\alpha(x,y)}} \Psi(x,y)
\left\{1+O\left(\frac{1}{u} +\frac{\log y}y\right)\right\},\ {\rm when\ }
1\ \leq\ d\ \leq\ y\ \leq\ {x \over d}.
\end{equation}
\bigskip

\begin{proposition}\label{Proposition_2.1} Throughout the range
(\ref{2.2}), for any $1 \leq d \leq x$, we have
$$
 \Psi\left( \frac xd,y\right) \leq \frac 1{d^{\alpha(x,y)}} \Psi(x,y)
\{ 1+o(1)\} ,
$$
where $\alpha$ is the solution to (\ref{alpha_def}).
In fact,
$$
\Psi \left ( {x \over d},y\right )\ <\  {\Psi(x,y) \over d^{\alpha(x.y)}},
$$
provided that\\
$$
\frac{\log d}{\log u\log y + \sqrt{ u \log u\log y}} \to \infty \, .
$$
\end{proposition}
\bigskip

\noindent {\bf Proof}.
By (\ref{2.1}), for $d=y^r$ with $0 \leq r \leq u/2$,  we have
$$
{\Psi\left( \frac xd,y\right)d^{\alpha} \over  \Psi(x,y)}\ =\
{d^{-(1-\alpha)} \rho(u-r) \over \rho(u)} \
\left(1+O\left({\log(u+1)\over \log y}\right)\right).
$$
The logarithm of the main term on the right side is
$$
-(1-\alpha)r \log y + \log (\rho(u-r)/\rho(u)).
$$
Using the fact that $u = (\log x)/(\log y)$, this can be rewritten
as
$$
r(\xi(u) - (1-\alpha)\log y) + \left( -\int_{u-r}^u \frac{\rho'(v)}{\rho(v)}dv
 -r\xi(u)\right) .
$$
The first term is $O(r/u)$ by (\ref{2.5}). Corollary 8.3 of \cite{tenen} gives that
\begin{equation} \label{2.8}
-\rho'(v)/\rho(v) = \xi(v) (1+O(1/v))\,,
\end{equation}
so that the second term equals
$$
 -\int_0^r (\xi(u)-\xi(u-t)) dt +O(r\log u/u).
$$

Now, differentiating $e^\xi=u\xi+1$ we obtain
$$
\xi+u\xi'\ =\ \xi' e^\xi\ =\ \xi' (u\xi+1),
$$
so that
$$
\xi'\ =\ {1 \over u - (u-1)\xi^{-1}}\ =\ {1 \over u(1 + O(1/\log
u))} =\ {1 \over u}   \left(1+O\left(\frac 1{\log u}\right)\right)
.
$$
Therefore
\begin{eqnarray} \label{xi_integral}
\int_0^r (\xi(u)-\xi(u-t)) dt &=& \int_0^r (r-v)\xi'(u-v) dv
= \left(1+O\left(\frac 1{\log u}\right)\right)\int_0^r \frac{(r-v)}{(u-v)} dv
\nonumber \\
&=& \left(1+O\left(\frac 1{\log u}\right)\right) (r-(r-u)\log (1-r/u))\,.
\end{eqnarray}

\noindent
Combining this with the above
yields that
\begin{eqnarray}
\log\left( {\Psi\left( \frac xd,y\right)d^{\alpha} \over \Psi(x,y)} \right)\ &=&\
- \left ( 1 + O \left ( {1 \over \log u} \right ) \right ) (r - (r-u) \log (1 - r/u))\nonumber \\
&&\hskip1.5in +\ O \left ( {r \log u \over u} +\frac{\log
(u+1)}{\log y} \right)
\nonumber \\
&=&\ -\frac{r^2}{2u} \left\{ 1 + O\left( {r \over u} + \frac
1{\log u} +\frac {\log u}r \right) \right\} + O\left( \frac{\log
(u+1)}{\log y} \right). \nonumber
\end{eqnarray}
 From (\ref{xi_integral}) and the first equation here we find
that this is negative provided
$r \leq u/2$ and $(\log u + \sqrt{u \log u / \log y}) / r \to 0$,
and is $o(1)$ in the complementary range.

If $d>\sqrt{x}$ we simply iterate the above result: The proposition follows
by noting that $\alpha(x,y)$ is a decreasing function in $x$ for fixed $y$,
by definition. \hfill $\Box$
\bigskip

We will require the following lemma, which is in one sense stronger, and
in another sense weaker, than Lemma \ref{first_lemma}.

\begin{lemma} \label{Lemma_2.2}  We have
$$
{\Psi(x,y) \over y}\ = o\left( {\Psi(x,y_0) \over y_0(\log y_0)^{1+\epsilon/4}
}  \right)
$$
for all $y$ outside of  the range
\begin{equation}
y_0 \exp( -(1+\epsilon) \sqrt{\log y_0\log\log y_0})\ \leq\ y\ \leq\
y_0 \exp( (1+\epsilon) \sqrt{\log y_0\log\log y_0});
\end{equation}
and
$$
{\Psi(x,y) \over y}\ \leq\ {(2/e^2-\epsilon) \Psi(x,y_0) \over y_0}
$$
for all $y$ outside of  the range
\begin{equation}
y_0 \exp( -(c_3+\epsilon) \sqrt{\log y_0})\ \leq\ y\ \leq\
y_0 \exp( (c_3+\epsilon) \sqrt{\log y_0}).
\end{equation}
\end{lemma}

\noindent {\bf Proof}.  Let $x=y_0^{u_0}$.
Define $g(u)=g_x(u)=\log \rho(u) - u^{-1} \log x$. By (\ref{2.1}) we have
$\log (\Psi(x,y)/xy) = g(u) +O(1/u)$, provided $\log y\asymp \log L$.
Select $u_1$
to maximize $g(u)$. Therefore $g(u_1)\geq g(u_0)$ by definition of $u_1$;
and $g(u_0)\geq g(u_1)+O(1/u_0)$ by the definition of $u_0$ and the above
estimate; therefore $g(u_0)=g(u_1)+O(1/u_0)$.

\

\noindent
By (\ref{2.8}), we have
$g'(v) = \rho'(v)/\rho(v) +v^{-2}\log x=-\xi(v)+v^{-2}\log x +O(\log v/v)$;
so that, for $t=O(u_1/\log u_1)$,
\begin{eqnarray}
g'(u_1+t)&=& g'(u_1+t)-g'(u_1)\nonumber \\
&=& \xi(u_1)-\xi(u_1+t)+\left(\frac1{(u_1+t)^2}-\frac1{u_1^2}\right)\log x
+O\left( \frac{\log u_1}{u_1}\right)\nonumber \\
&=& O\left( \frac{t+\log u_1}{u_1}\right)
-2tu_1^{-3}\log x
(1+O(t/u_1))  \nonumber \\
&=& -2t\frac{\xi(u_1)}{u_1} + O\left( \frac{t+\log u_1}{u_1}\right), \nonumber
\end{eqnarray}
since
$0=g'(u_1)= -\xi(u_1) +u_1^{-2}\log x+O(\log u_1/u_1)$. Therefore
\begin{equation} \label{eq:gu1}
g(u_1)-g(u_1+T) = - \int_0^T g'(u_1+t) dt = \frac{T^2}{u_1} (\xi(u_1)+O(1))
+O\left( \frac{T\log u_1}{u_1}\right) ,
\end{equation}
for $T=O(u_1/\log u_1)$. We deduce that
$u_0=u_1+O(1)$, as well as
both
$$
g(u)<g(u_0)-(1+\epsilon/3) \log u_0 \mbox{ for \ } |u-u_0|>(1+\epsilon/2) \sqrt{u_0}\,,
$$
and
$$
g(u)<g(u_0)-\log(e^2/2+\epsilon) \mbox{ for \ } |u-u_0|>(c_3+\epsilon)
\sqrt{u_0/\log u_0}\, ,
$$
which are the  desired
results. \hfill $\Box$

\bigskip

Next we obtain a more accurate estimate for $y_0$ than (\ref{J0y0}):

\begin{lemma} \label{Lemma_2.3} We have
\begin{eqnarray*}
\log y_0 & = & \log L(x) \left( 1 + \frac{\log_3x-\log 2}{2\log_2x} 
   + O\left( \left( \frac{\log_3x}{\log_2x}\right)^2 \right)\right)  
   \;\; \mbox{ and } \\[1ex]
\frac{u_0 \xi(u_0)}{\log y_0} & = & 
   1+ O\left( \frac 1{u_0}\right) \, . 
\end{eqnarray*}
\end{lemma}

\noindent {\bf Proof}. In the notation of the Lemma \ref{Lemma_2.2} we see 
by (\ref{eq:gu1}) that $|g(u_1+T)| = o(1/u_1)$ as $T \to \infty$, 
so that $u_0=u_1+O(1)$.
We saw that $u_1^2 \xi(u_1) (1+O(1/u_1)) = \log x$, so the same equation is satisfied by
$u_0$ (in place of $u_1$), and the estimate for $\log y_0=(1/u_0)\log x$ follows
from (\ref{2.4b}). Moreover $u_0 \xi(u_0)= \log y_0 (1+O(1/u_0))$
\hfill $\Box$
\bigskip

\begin{corollary}\label{Corollary_2.4} If $d=p_1p_2\ldots p_k$, where
each $p_j$ is a prime in $(y,My]$ we have
\begin{equation} \label{eq:unif}
\frac{\psi (x / (p_1 \cdots p_k) \, , \, y_0)}{\psi (x , y_0)}
   \sim \frac{(\log y_0)^k}{p_1 \cdots p_k}
\end{equation}
uniformly in $k \geq 1$ and $\log M=o((\log x/\log\log x)^{1/4})$,
as $x \to \infty$.  Also
\begin{equation} \label{eq:unif unif}
\frac{\psi (x / (p_1 \cdots p_k) \, , \, y_0)}{\psi (x , y_0)}
   \leq 2^k \frac{(\log y_0)^k}{p_1 \cdots p_k}
\end{equation}
uniformly for  $k \geq 1$ and $\log M = o(\log x/\log\log x)^{1/2}$, as
$x \to \infty$.
\end{corollary}

\noindent {\bf Proof}.  We use (\ref{2.7}) at most $2k$ times to obtain
$$
\frac{\psi (x / (p_1 \cdots p_k) \, , \, y_0)}{\psi (x , y_0)} =
\frac 1{(p_1 \cdots p_k)^{\alpha}}
\left\{1+O\left(\frac{k}{u_0} +\frac{k\log y_0}{y_0}\right)\right\}
\sim \frac  {(p_1 \cdots p_k)^{1-\beta}} {p_1 \cdots p_k}
$$
where $\alpha(x,y_0)\geq \beta\geq \alpha(x/(p_1 \cdots p_k),y_0)$.
If $u'=\log(x/(p_1 \cdots p_k)/\log(y_0)$ then $u'=u+O(k)$ and so
$y_0^{1-\beta}=u_0\xi(u_0)\{ 1+O(k/u_0)\}=\log y_0\{ 1+O(k/u_0)\}$, by (\ref{2.6})
and then Lemma \ref{Lemma_2.3}. Hence we obtain (\ref{eq:unif}) as $k^2=o(u_0)$
and, in our range,
$$
M^{k(1-\alpha)}=\exp( O(k \log M (\log\log y_0)/( \log y_0)))=1+o(1).
$$
To obtain (\ref{eq:unif unif}) we can use the same estimates but now we simply need
$k / u_0 \to 0$ so that $y_0^{1-\beta}\leq (4/3) \log y_0$, and
$\log M  / u_0$ so that $M^{1-\beta}\leq (4/3)$.
\hfill $\Box$

\subsection{Straightforward analytic  estimates}

We complete this section by collecting together various straightforward
analytic  estimates that will be needed later.

Fix $0 < a < b$.
By the prime number theorem, we  have
\begin{equation} \label{eq:PNT 1}
\sum_{ay < q \leq by } \frac{\log y}{q}
   \sim \log \left ( \frac{b}{a} \right ) \, .
\end{equation}
where the sum is over primes $q$, and also that
\begin{equation} \label{eq:PNT 2}
\sum_{ay < q \leq by} \frac{\log y}{q}
   \leq 2 \log \left ( \frac{b}{a} \right ) \, ,
\end{equation}
for all $1 \leq a \leq b/2$, once $y$ is sufficiently large.
To see this note that, since \\
$\sum_{q\leq Q} (\log q)/q=\log Q+C+o(1)$,
for some constant $C$,  the sum is
$$
\leq \sum_{ay < q \leq by} \frac{\log q}{q} =\log \left ( \frac{b}{a} \right ) +o_{y\to \infty}(1),
$$
and the result follows.

\medskip

\begin{lemma} \label{Lemma_integral1}  Let
\begin{equation} \label{eq:def g}
  g(\beta,C)\ :=\   \beta^{-2} \int_0^{C/\beta^2} \log\left(
\frac{e^{z} + e^{-z}}2 \right)  {dz \over z^2} + 1- \log(C)
.
\end{equation}
The function $g(1,C)$ is decreasing for $C>0$, with
$$
\lim_{ C\to \infty}  g(1,C)  = \gamma+\log(4/\pi)\,.
$$
\end{lemma}

\noindent {\bf Proof}. Since
$$
\frac{dg(1,C)}{dC}    =
 {\log( \frac 12 (e^C + e^{-C})) \over C^2}\ -\ {1 \over C}\ < 0,
$$
for all $C>0$,  we minimize by letting $C\to \infty$.  Integrating by parts,
we have that
$$
\lim_{ C\to \infty}   g(1,C)  =
\int_0^1
\frac{e^{z}-e^{-z}}{e^{z}+e^{-z}} \frac{dz}{z }
-2\int_1^\infty \frac{e^{-z}}{e^{z}+e^{-z}} \frac{dz}{z} .
$$
Now 6.1.50 of \cite{AS} states that
$$
\log \Gamma(s)= \int_0^\infty \left( (s-1) e^{-t} - \frac{e^{-t}-e^{-st}}
{1-e^{-t}}\right)  \frac {dt}t;
$$
and  the third line of 6.3.22 of \cite{AS} readily implies that
\begin{equation} \label{gamma}
\gamma =  \int_0^1   (1-e^{-t}) \frac{dt}t
 -\int_1^\infty  e^{-t} \frac{dt}t
.
\end{equation}
Since $\Gamma(1/2)=\pi^{1/2}$, and taking $s=1/2$ and $t=4z$, our result follows. \hfill $\Box$

\section{The lower bound for $T$ in Theorem \ref{main_theorem}, and
Theorem \ref{Theorem_2} } \label{section5}

\subsection{Proof strategy}

To establish that
$$
\PP\Bigl(T\ >\ (\pi/4)(e^{-\gamma} - \epsilon) J_0(x)\Bigr)\ =\ 1 - o(1),
$$
we show  that the expected number of non-trivial subsets $S$ of
$\{ 1,...,J\}$ for which $\prod_{i\in S} a_i$ is a square is
$o(1)$, for $J(x) = (\pi/4)(e^{-\gamma}-o(1)) J_0(x)$.

\subsection{Structure of a square product}

We begin with the following proposition.

\begin{proposition} \label{Proposition_5.1}
Select integers $a_1,\dots , a_J$ at random from $[1,x]$.
The probability that there exists a subsequence $I$ of the $a_i$
with
$$
2\ \leq\ |I|\ \leq\  {\log x \over 2\log\log x} \ \textit{  for\
which\ } \prod_{a \in I} a\  \textit{is a square}
$$
is $O(J^2\log x/x)$ provided $J< x^{o(1)}$.
\end{proposition}

\noindent {\bf Proof}.  Suppose that $b_1,\dots ,b_k$ were chosen at random  from $[1,x]$.
The probability that  $b_1b_2\dots b_k$ is a square equals
$$
x^{-k} |\{ b_1,\dots ,b_k\ \leq\  x\ :\ b_1b_2\dots b_k \in \mathbb Z^2\}|.
$$
Now write each $b_i$ uniquely as
$$
b_i\ =\ c_iu_i^2,\ {\rm where\ }c_i\ {\rm is\ squarefree}.
$$

Assuming that $b_1 \cdots b_k$ is a square, which implies
$c_1 \cdots c_k$ is a square, define the doubly indexed sequence
$c_{i,j}$, where $i,j=1,...,k$ and $i \neq j$, to be any satisfying the
relations
\begin{equation} \label{ci_properties}
c_{i,j}\ =\ c_{j,i}, \ \ {\rm with\ }\ c_i\ =\ \prod_{j \neq i}
c_{i,j}  \ \ {\rm for\ each \ } i.
\end{equation}
The fact that such $c_{i,j}$ exist can be seen as follows:  For each prime $p$
dividing $c_1 \cdots c_k$, we will need to decide which $c_{i,j}$ that $p$
divides; and, to do this, suppose that $p$ divides $c_{i_1},...,c_{i_{2t}}$
(the reason it is $2t$ is that all the $c_i$ are square-free and have product
a square).  Then, the following $c_{i,j}$ are to be divisible by $p$, and no others:
$$
c_{i_1,i_2},\ c_{i_2,i_1},\ c_{i_3,i_4},\ c_{i_4,i_3},\ ...,\ c_{i_{2t-1}, i_{2t}},\ c_{i_{2t},
i_{2t-1}}.
$$
Each $c_{i,j}$ is then the product of the primes dividing $c_1
\cdots c_k$ which divide it; and if this process leaves some
$c_{i,j}$ not divisible by any prime $p | c_1 \cdots c_k$, then we
set $c_{i,j} = 1$.
\bigskip

Given $c_1,...,c_k$, the number of sequences $b_1,...,b_k$ satisfying
$b_i = c_i u_i^2$ is the
number of possibilities for the numbers $u_i$, which is
$\leq (x/c_i)^{1/2}$; and so, the probability that  $b_1\cdots b_k$ is a square is
\begin{eqnarray}
&\leq& \frac 1{x^k} \sum_{c_{i,j}\leq x \atop {\rm for} \ 1\leq i<j\leq k}
\prod_{i=1}^k \left(  \frac{x} {\prod_{j\ne i} c_{i,j}} \right)^{1/2} \nonumber \\
&\leq& \frac 1{x^{k/2}} \sum_{1\leq i<j\leq k} \left( \sum_{c_{i,j}\leq x}
  \frac{1} {c_{i,j}} \right)  \leq \frac 1{x^{k/2}}\ (1 +\log x)^{k(k-1)/2}
\end{eqnarray}
since each $c_{i,j}$ appears twice in the above product. Therefore the probability
that there exists $I\subset \{ 1,2,\dots ,J\}$  for which
$\prod_{i\in I} a_i \in \mathbb Z^2$ with $|I|=k$ is
$$
\leq {J \choose k} \frac 1{x^{k/2}}\ (1 +\log x)^{k(k-1)/2}
\leq \left( \frac{ J^2 (1 +\log x)^{k-1}}{x} \right)^{k/2}
$$
which gives $O(J^2\log x/x)$ for $k=2$, and is
$\leq 1/x$ for $3\leq k\leq \log x/2\log\log x$. \hfill $\Box$

\subsection{The main argument}

In this subsection, we prove that
$$
\PP\Bigl(T\ >\ (\pi/4)(e^{-\gamma} - \epsilon) J_0(x)\Bigr)\ =\ 1-o(1).
$$

As a consequence of the upper bound proved in \cite{CGPT}, we may assume that
$T < (3/4)J_0(x)$ holds with probability $1-  o(1)$. Furthermore, following Proposition
\ref{Proposition_5.1} we need only focus on
subsequences $I$ of $a_1,...,a_J$ (where $J = T < J_0(x)$) of
length exceeding $\log x/2 \log\log x$, that have product equal to
a square.
\bigskip

Throughout we shall write $a_i=b_id_i$ where $P(b_i)\leq y$ and where either
$d_i = 1$ or $p(d_i) > y$ , for
$1\leq i\leq k$.  Recall here that $p(n)$ denotes the smallest and $P(n)$ the largest prime divisor of $n$. If $a_1,\dots ,a_k$ are chosen at random from $[1,x]$
then
\begin{eqnarray} \label{5.1}
\PP(a_1\dots a_k \in \mathbb Z^2)\ &\leq&\ \PP(d_1\dots d_k \in \mathbb Z^2)\ \nonumber \\
&=&\ \sum_{d_1,\dots, d_k\geq 1 \atop {d_1\dots d_k \in \mathbb Z^2 \atop d_i=1\ {\rm or\ }p(d_i)>y}}
\prod_{i=1}^k \frac{\Psi\left( x/d_i,y\right)}x \nonumber \\
&\leq& \left( \{ 1 +o(1)\} \frac{\Psi(x,y)}x \right)^k
\sum_{n = 1\ {\rm or\ }p(n)>y} \frac{\tau_k(n^2)}{n^{2\alpha}}\,,
\end{eqnarray}
by Proposition \ref{Proposition_2.1}, where $\tau_k(m)$ denotes the number
of different ways of writing $m$ as the product of $k$ positive integers.

Out of $J=\eta J_0$ integers, the number of $k$-tuples is ${J \choose k}\leq (eJ/k)^k$;
and so the expected number of $k$-tuples whose product is a square is
\begin{equation}\label{5.2}
\leq \left(  (e+o(1)) \frac{ \eta y}{k\log y_0} \frac{\Psi(x,y)/y}{\Psi(x,y_0)/y_0}
\right)^k \prod_{p>y} \left( 1 + \frac{\tau_k(p^2)}{p^{2\alpha}}
+ \frac{\tau_k(p^4)}{p^{4\alpha}} + \dots \right) \,.
\end{equation}
\bigskip

We now consider $k$ in two different ranges, and in both ranges we
will select different values for $y$, so as to give good upper bounds for
(\ref{5.2}):
\bigskip

$\bullet$  First, if
$$
{\log x \over 2\log\log x}\ <\ k\ \leq\ y_0^{1/4},
$$
then let $y=y_0^{1/3}$ so that
$k=o(y_0^\alpha)$. Therefore the Euler product in  (\ref{5.2}) is
$$
\leq \exp \left( O\left( \sum_{p>y}\frac{k^2}{p^{2\alpha}}  \right) \right)
\leq \exp \left( O\left(  \frac{k^2y^{2(1-\alpha)}}{y \log y}  \right) \right)
=e^{o(k)}.
$$
Now $\Psi(x,y_0^\gamma) = x/y_0^{1/\gamma+o(1)}$ by
(\ref{psi_estimate_2}) and therefore the quantity in (\ref{5.2})
is
\begin{equation}\label{goodbound}
\leq \left( \frac{1/y_0^{3+o(1)}} {k/y_0^{2+o(1)}}\right)^k \leq
y_0^{-k+o(k)} ,
\end{equation}
which is $<1/x^2$ in this first
range for $k$.
\bigskip

$\bullet$ Next, we consider the range
$$
y_0^{1/4}\ \leq\ k\ =\ y_0^\beta\ \leq\ J\  \leq\ J_0.
$$
In this case we will choose
$y$ so that $[k/C] = \pi(y)$, and then will optimize the $C$ later.  For this choice of
$y$ a simple calculation reveals that
\begin{eqnarray}
1 + {\tau_k(p^2) \over p^{2\alpha}} + {\tau_k(p^4) \over p^{4\alpha}} + \cdots
\ &\sim&\ 1 + {(k/p^\alpha)^2 \over 2!} + {(k/p^\alpha)^4 \over 4!}  + \cdots \nonumber \\
&=&\ {e^{k/p^\alpha} + e^{-k/p^\alpha} \over 2}. \nonumber
\end{eqnarray}
In order to evaluate (\ref{5.2}) we need to product this over primes
$p > y$.  The logarithm of this product equals
$$
\sum_{p > y \atop p\ {\rm prime}} \log\left( \frac{e^{k/p^\alpha}
+ e^{-k/p^\alpha}}2 \right) \ \sim\ \int_y^\infty {1 \over \log t}
\log\left( \frac{e^{k/t^\alpha} + e^{-k/t^\alpha}}2 \right)  dt,
$$
by the prime number theorem. Letting $z = k/t^\alpha$, from
(\ref{y_alpha}) this last integral is
$$
\sim\ \int_0^{C/\beta^2}
 {(k/z)^{1/\alpha} \over z\log(k/z)} \ \log\left( \frac{e^{z} + e^{-z}}2 \right) dz.
$$
Now, $k^{1/\alpha} \sim \beta^{-2}\log y$ by (\ref{y_alpha}) so
that
$$
\frac{(k/z)^{1/\alpha}}{\log (k/z)}\  \sim\ (k/z) \beta^{-2}
$$
as $z = o(1)$.  It follows that the quantity in (\ref{5.2}) 
is bounded from above by
\begin{equation} \label{gbC}
\left ( (1 + o(1))  e^{ g(\beta,C)} \beta \eta
{\Psi(x,y)/y \over \Psi(x,y_0)/y_0} \right )^k,
\end{equation}
where $g(\beta,C)$ is defined in (\ref{eq:def g}).
\bigskip

Now, for any fixed $C$ we have, as a consequence of Lemma~\ref{first_lemma}, 
that (\ref{gbC}) is $o(1/x^{2})$ unless $\beta =
1 + o(1)$; and so, we really only need to consider $k = y_0^{1 +
o(1)}$, as the total expected number of $k$-tuples for other
values of $k$ add only  $o(1/x^{2+o(1)})$.  If $C=C(\epsilon)$ is sufficiently
large then $ e^{g(1,C)} <4e^\gamma/\pi +\epsilon$ by
Lemma~\ref{Lemma_integral1} and,  since $y_0$ maximizes
$\Psi(x,y)/y$ for $y=y_0$, we deduce that (\ref{5.2}) is at most
$$
\leq ( (1+\epsilon) 4\eta  e^\gamma/\pi)^k.
$$
Therefore, if $\eta< (1-\epsilon)  e^{-\gamma} \pi/4$, then this is
less than $1/x^2$. \hfill $\Box$

\subsection{Proof of Theorem~\ref{Theorem_2}, part  (a)}

This last proof yields further useful information: If
either $J<(\pi/4)(e^{-\gamma} - \epsilon) J_0(x)$, or if $k<y_0^{1 -o(1)}$ or $k>y_0^{1 +o(1)}$, then the expected number of square products with $k>1$ is
$O(J_0(x)^2\log x/x)$, whereas the expected number of squares in our sequence
is $\sim J/\sqrt{x}$. This justifies  Theorem~\ref{Theorem_2}(a).

\subsection{Proof of Theorem \ref{Theorem_2}, part  (b)}

The proof in Section~3.3 yielded that if we have a square product then,
with probability $1+o(1)$,  we have $|I|=k=y_0^{1+o(1)}$.
We now assume that $k=y_0^{1+o(1)}$ with
\begin{equation} \label{k_range}
 k\ \not\in\
[ y_0 \exp(-(c_3 + \epsilon)\sqrt{\log y_0}),\ y_0 \exp( (c_3 + \epsilon) \sqrt{\log y_0})].
\end{equation}
>From the discussion following  (\ref{gbC}) above,   we know,  by taking $C$ large,
that the number of such $k$-tuples is at most
$$
\left ( (4e^\gamma/\pi+\epsilon) {\eta \Psi(x,y)/y \over
\Psi(x,y_0)/y_0} \right )^k.
$$
By Lemma \ref{Lemma_2.2}, this is at most
$$
\left ( (4e^\gamma/\pi + \epsilon)(2/e^2 + o(1)) \eta \right )^k\
<\ 1/2^k,
$$
for sufficiently small $\epsilon > 0$, using the fact that $\eta < 3/4$.
Therefore the expected number of $k$- tuples with product a square is
$o(1)$ for all $k$ satisfying (\ref{k_range}), so that
Theorem~\ref{Theorem_2}(b) follows. \hfill $\Box$

\subsection{Proof of Theorem \ref{Theorem_2}, part (c)}

In the previous subsection we proved that
$$
|I|\ \leq\  y_1\ :=\ y_0 \exp(  (1+\epsilon) \sqrt{\log y_0\log\log y_0}),
$$
with probability $1-o(1)$. In this section we prove, among other
results,  part (c) of Theorem~\ref{Theorem_2}.

\begin{proposition} \label{Proposition_5.2}
Write each $a_i=b_id_i$ where  $P(b_i)\leq y=y_1<p(d_i)$, and
suppose that $d_{i_1}\dots d_{i_l}$ is a subproduct which equals a
square $n^2$, but such that no subproduct of this is a square. Then,
with probability $1-o(1)$, we have $l=o(\log y_0)$ and
$n$ is a squarefree integer composed of precisely $l-1$ prime factors,
each $< y^2$, where $n\leq y^{2l}$.
\end{proposition}

\noindent {\bf Proof}. For ease of notation we will relabel,
replacing $d_{i_1}\dots d_{i_l}$ by $d_1\dots d_l$. Note that with
the choice of $y=y_1$, we have $y/l\log y\to \infty$ and
$y=y_0^{1+o(1)}$, so we know that $y^{\alpha}\sim y/\log y$ by
(\ref{y_alpha}).

We now show that $n$ has at least $l-1$ (not necessarily distinct)
prime factors, so that $n^2=d_{1}\dots d_{l} > y^{2(l-1)}$: Create a graph $G$
on the $l$ vertices $v_1,\ldots, v_l$ where, for each prime $p^q$
which (exactly) divides  $n$, draw a total of $q$ edges, placing an edge 
between pairs of vertices $v_j$ for which $p$ divides $d_j$. Now $G$ is
connected, since our square product is minimal, and so must have $\geq
l-1$ edges.

We now modify the argument from the start of section~3.3 (with $k$
replaced by $l$) to restrict our attention to cases in which
$d_1\dots d_l\geq y^{2l}\phi(x)^2$, where $\phi(x)=y^{O(1)}$. 
\noindent
To obtain an upper bound we may multiply through the summand, in
(\ref{5.1}), by $(n/y^{l}\phi(x))^{2\theta}$, where we have chosen
$\theta>0$ so that $y^{2\theta}=(2y\log l)/(l(\log y)^2)$. Then we 
must multiply the right side of (\ref{5.2}) through by $1/(y^{2\theta})^l
\phi(x)^{2\theta}$ and change the terms in the Euler product to
$(1 + \tau_l(p^2)/p^{2(\alpha-\theta)} + \tau_l(p^4)/p^{4(\alpha-\theta)} 
+ \dots )$.  


First we bound the Euler product using the prime number theorem:
Recall that  the function $\tau_\ell(n)$ 
counts the number of sequences of positive integers $d_1,...,d_\ell$
such that
$
d_1 \cdots d_\ell\ =\ n.
$
In the case $n = p^{2k}$, this amounts to computing the number of ordered
partitions of $2k$ into $\ell$ parts that are $\geq 0$; so, 
$$
\tau_\ell(p^{2k})\ =\ {2k + \ell - 1 \choose 2k}\ \leq\ 
\left \{ \begin{array}{rl} \ell(\ell+1)/2,\ &{\rm if\ } k=1; \\
(2k+\ell-1)^{2k} \over (2k)!,\ &{\rm if\ } k \geq 2.\end{array}\right.
$$
For
$
p\ =\ y_0^{1+o(1)}\ =\ L(x)^{1+o(1)}$,   using (\ref{y_alpha}) with $\beta=1$, we have that

$$
{1 \over p^{2\alpha}}\ \sim\ {\log^2 p \over p^2} 
$$
making the summation of terms involving $p$ in the Euler product become: 
$$
\{ 1+o(1)\}  \frac{\ell (\ell+1)}2 \cdot \frac{\log^2p}{p^2} \cdot 
p^{2\theta}\,.
$$
Via the prime number theorem the logarithm of the Euler product is therefore
$$
\sim \frac{\ell (\ell+1)}2  \sum_{y<p<y^4} \frac{\log^2p}{p^{2-2\theta}} 
\sim \frac{\ell (\ell+1)}2  \int_{y}^{y^4} \frac{\log t}{t^{2-2\theta}}  dt.
$$
(Here the  primes $p$,
with $y < p < y^{4+o(1)}$,  being the only relevant ones follows 
from  comments made above the statement of Theorem~\ref{Theorem_2}.)
Now $\theta\leq 1/2$ by definition, so the above calculation becomes
$$
\sim \frac{\ell (\ell+1)}2  \cdot \frac{\log y}{(1-2\theta) y^{1-2\theta}}   
= \frac{\ell (\ell+1)}2  \cdot \frac{\log^2 y}{y^{1-2\theta} (1-2\theta)\log y }\,.
$$
Now $y^{1-2\theta}=\ell \log^2 y/2\log \ell$, so the above is
$$
=   \frac{ (\ell+1)   \log \ell }
{  \log \left( \ell \log^2 y/2\log \ell\right) }
= \ell \left( 1  + O\left( \frac{\log\log \ell}{\log \ell}\right) \right) 
= \{ 1+o_{\ell\to\infty}(1)\} \ell\,.
$$
So putting (\ref{5.2}) to use as explained above,  the  expected number of such $l$-tuples is
\begin{equation} \label{5.39}
\leq \frac 1{\phi(x)^{2\theta}} \left(  (e+o(1)) \ \frac{ \eta y}{\ell y^{2\theta}\log y_0}\
\frac{\Psi(x,y)/y}{\Psi(x,y_0)/y_0} \right)^l \ e^{(1+o(1)) \ell} 
\end{equation}

\begin{equation} \label{5.399}
=  \frac 1{\phi(x)^{2\theta}} \left(   (e+o(1)) \ \frac{\eta (\log^2 y)}{(2\log \ell)(\log y_0)}\
\frac{\Psi(x,y)/y}{\Psi(x,y_0)/y_0} \right)^l \ e^{(1+o(1)) \ell} 
\end{equation}

\begin{equation} \label{5.4}
 \ll
\frac {1}{\phi(x)^{2\theta}(\log y_0)^{\epsilon l/5}}  \,,
\end{equation}
as $\eta \le 1$, and by Lemma~{\ref{Lemma_2.2}} for $y=y_1$.

Now we are ready to establish the conclusions of the proposition.
Take $\phi(x)=1/y$ in the above, and as $2\theta<1$ by definition,
(\ref{5.4}) becomes $\ll y /(\log y)^{\epsilon \ell/5}$\,.
This is $o(1)$ provided $\ell \ge 6 \log y/(\epsilon \log\log y) $,
hence we expect $o(1)$ products with $l\gg \log y_0$,
yielding $l=o(\log y_0)$ with probability
 $1-o(1)$. In this case $2\theta\sim 1$. 

Regarding the structure of the factorization of $n$:
Taking $\phi(x)=1$, we  expect  $o(1)$ products with $d_1\dots d_l\geq y^{2l}$; hence
$d_1\dots d_l=n^2<y^{2l}$ with probability
$1-o(1)$. Since each prime divisor is $>y$, evidently $n$ has $<l$ prime factors,
and so exactly $l-1$. Also,  if $p$ is the largest then  $y^{l-2}p<y^{l}$,
that is $p<y^2$.

Finally, we are left with showing that $n$ is squarefree.
To obtain an upper bound on
the expected number of square products $n^2$ for which $n$ is divisible
by the square of a prime $>y$, we proceed  much as above
with $\phi(x)=1/y$, but now
the Euler product has an additional factor
$$
\sum_{p>y} \left(  \frac{\tau_l(p^4)}{p^{4\alpha}} +
\frac{\tau_l(p^6)}{p^{6\alpha}} + \dots \right)  \ll
\frac{l^4}{(y/\log y)^3}
 \ll \frac {(\log y)^7}{y^3} .
$$
 From  (\ref{5.4}) we thus deduce that we expect $o(1)$ such square products.

$\Cox$

\section{Hypergraphs} \label{sect:Hyper}

The main result of this section is to prove the upper bound in 
Theorem~\ref{main_theorem}.  A roadmap for the proof is as
follows.  

Recall that the numbers $a_1 , a_2 , \ldots$, chosen uniformly at random from $\{1,2, \ldots, x\}$,  are encoded as row 
vectors over $\F_2$.  Subsets whose product is a square are determined
by combinatorial relations among these row vectors.  Schroeppel's
method and its variants ignore columns corresponding to primes less
than $y_0$.  This makes the relations easier to satisfy but we pay
for it by requiring $\pi (y_0)$ many relations.  To make the search
more tractable, we restrict our attention to the more obvious ways
of finding linear relations.  Schroeppel's original method considers
only the most obvious: after removing columns less than $y_0$ we must
be left with all zeros.  The {\sl one large prime} variation considers
also the next most obvious: when we have two identical rows containing
a single~1.  

The upper bound in Theorem~\ref{main_theorem} is proved via the 
{\sl $k$ large primes variation}.  
We consider only rows in which at most $k$ ones remain.
Tractability of the analysis rests on the fact that the combinatorial
structure converges as $x \to \infty$ to a random object built from a 
Poisson point process.  In order for the convergence to  be uniform, in 
addition to restricting $k$, we must restrict the columns: specifically,
fixing $M > 0$, we must not use any $a_i$ with a prime factor 
greater than $M y_0$.  We must also restrict the combinatorial
complexity of the search for linear relations as follows: calling
two rows ``neighbors'' if they share a nonzero column (whose index
is now forced to be between $y_0$ and $M y_0$), any linear relation
must take place within a ball of some fixed radius $m$ in the 
neighbor graph on rows.  We may then prove that the combinatorial
structure converges in an appropriate sense to a tree-like random
hypergraph defined on a Poisson point process.  The number of samples
needed to accumulate $\pi (y_0)$ linear relations in the limiting
model is computable explicitly in terms of some functions $\gamma_{m,M,k}$.
For fixed $m,M,k$, these are ugly, but as $m,M,k \to \infty$, this
number decreases to $e^{-\gamma} J_0$.

An outline of this section is as follows.  Section~\ref{ss:4.1}
defines some functions that include the family $\{ \gamma_{m,M,k} \}$.  
A result (Theorem~\ref{th:quant}) is then formulated in terms of these 
functions which implies the upper bound in Theorem~\ref{main_theorem}.
The subsection ends with the definition of some combinatorial structures
such as tree-like hypergraphs that will be used in the search for 
linear relations.  Section~\ref{ss:4.2} formally defines the probability
model and the random objects (hypergraphs with distinguished vertices)
that will witness linear relations.  The number of rows 
neighboring any given row is shown to have finite first
and second moments (Proposition~\ref{pr:counting}), which is
then parlayed into an upper bound on the mean of size of the
$m$-ball in the neighbor graph on rows.  Section~\ref{ss:4.3}
constructs the limit object, an informal description of which
appears at the beginning of that subsection.  Section~\ref{ss:4.4}
proves convergence of the random hypergraphs in Section~\ref{ss:4.2}
to the limit object of Section~\ref{ss:4.3}.  Although it takes
several pages, it consists merely of repeated applications of
Proposition~\ref{pr:counting}.  Section~\ref{ss:4.5} evaluates
the probability $\theta_{m,k}^{M,\eta} (\rho)$, which is the 
probability in the limit model that if a row containing a single~1
in column $\rho y_0$ arises at time $\eta J_0$, it will form
a new linear relation.  The key result here (Lemma~\ref{lem:theta converge})
is that this is~1 when $m,M,k$ are sufficiently large and 
$\eta > e^{-\gamma}$.  Finally, Section~\ref{ss:4.6} finishes
the proof  of the main theorems.

\subsection{Preliminary results}
\label{ss:4.1}
 
To begin in earnest, we define the following functions, which will
arise in the branching processes with finite values of $m, k$ and $M$.
\begin{eqnarray*}
\exp_k (z) & := & \sum_{j=0}^{k-1} \frac{z^j}{j!} \, ;  \\
A_M (z) & := & \int_{1/M}^1 \frac{1 - e^{-zt}}{t} \, dt \, .
\end{eqnarray*}
Clearly, as $k , M \to \infty$, we have the limits
\begin{eqnarray*}
\exp_k (z) & \uparrow & \exp (z)  \, ;  \\
A_M (z) & \uparrow & A (z) := \int_0^1 \frac{1 - e^{-zt}}{t} \, dt \, .
\end{eqnarray*}
Recursively, define functions $\gamma_{m , M , k}$ for $m = 0, 1, 2, 
\ldots$ by
\begin{eqnarray}
\gamma_{0,M,k} (u) & := & u  \, ; \nonumber \\[1ex]
\gamma_{m+1,M,k} (u) & := & u \, \exp_k \left [ A_M (\gamma_{m,M,k}
   (u)) \right ] \, . \label{eq:gamma}
\end{eqnarray}
Note that $\gamma_{m,M,k} (u)$ is increasing in all four 
arguments.  From this it follows that
$\gamma_{m,M,k} (u)$ increases to $\gamma_{M,k} (u)$  as $m \to \infty$,
a fixed point of the
map $z \mapsto u \exp_k (A_M (z))$, so that
\begin{equation}\label{gamma_u}
\gamma_{M,k} (u)  := u \, \exp_k \left [ A_M (\gamma_{M,k} (u)) \right ] .
\end{equation}
We now establish that $\gamma_{M,k} (u)<\infty$ except perhaps 
when $M=k=\infty$: we have $0\leq A_M (z)\leq \log M$ for all $z$, 
so that $u<\gamma_{M,k} (u)\leq Mu$ for all $u$; in particular 
$\gamma_{M,k} (u)<\infty$ if $M<\infty$.  Also, $A(z)=\log z+O(1)$ 
which, along with~(\ref{eq:gamma}), implies that $\gamma_{\infty,k}(u) 
\sim u(\log u)^{k-1}/(k-1)!$; in particular $\gamma_{\infty,k} (u) <
\infty$.  As $M,k \to \infty$, the fixed point $\gamma_{M , k} (u)$ 
increases to the fixed point $\gamma (u)$ of the map 
$z \mapsto u e^{A(z)}$, or to $\infty$ if there is no such fixed
point, in which case we write $\gamma(u)=\infty$.  In
Lemma \ref{lem:theta converge} we show that this map has 
a fixed point if and only if $u \leq e^{-\gamma}$.  Otherwise 
$\gamma (u) = \infty$ for $u > e^{-\gamma}$ so that 
\begin{equation} \label{eq:eta star}
\int_0^\eta \frac{\gamma (u)}{u} \, du = \infty > 1
\end{equation} 
for any $\eta>e^{-\gamma}$.

\noindent
Our main result in this section is the following:

\begin{theorem} \label{th:quant}
If $\eta , m , M , k$ are such that
$$\int_0^\eta \frac{\gamma_{m+1,M,k} (u)}{u} \, du > 1, $$
then with probability approaching~1, as $x \to \infty$, among \,
$\eta J_0$ uniform random samples from $\{ 1 , \ldots , x \}$,
the $y$-smooth numbers up to $M y$ with at most $k$
large primes will contain a square subproduct.  Furthermore,
this will be witnessed in diameter at most $m$, in a sense
to be made precise in Definitions~\ref{def:chi-marked} and \ref{def:chi} below.
\end{theorem}
Together with~(\ref{eq:eta star}), this establishes the upper bound 
in Theorem~\ref{main_theorem}.  Our conjecture that the upper bound
is sharp is supported by the fact that $\lim_{t \uparrow \eta_*}
\int_0^t \frac{\gamma (u)}{u} \, du = 1$.

\subsubsection*{Hypergraphs}

A \Em{hypergraph} on a vertex
set $V$ is simply a collection $\hyp$ of finite subsets of $V$ of
cardinality at least~2.
Each $S \in \hyp$ is called a \Em{hyperedge} of $\hyp$; the
\Em{cardinality} of a hyperedge $S$ is its cardinality as a set.
Define the \Em{support} of a hypergraph $\hyp$, denoted
by $\supp (\hyp) := \bigcup_{S \in \hyp} S$, to be the
union of all of its hyperedges.  By a hypergraph $\hyp$
with vertex set $V$, we mean that $\supp (\hyp) \subseteq V$ (note:
in the literature, often this language would imply $\supp (\hyp) = V$).
We will typically use script letters for hypergraphs:
$\graph, \hyp$, and so forth.  A rooted hypergraph is simply a
hypergraph together with a choice of a distinguished element in
its support.  Thus, the hypergraphs on $V$ rooted at $p$ are in
one to one correspondence with hypergraphs on $V$ containing $p$
in their support.

\begin{defn}[tree-like hypergraphs]
A finite hypergraph $\graph$ rooted at $p$ is \Em{tree-like}
if $\supp (\graph)$ may be given the structure of a tree $T$,
rooted at $p$, in such a way that the following decomposition
holds.  Let $I$ denote the set of vertices that are not leaves
of $T$.  We require that for each $q \in I$,
the set of children of $q$ may be partitioned into sets $V_{q,1} ,
\ldots , V_{q , n(q)}$ so that each hyperedge of $\graph$
is equal to $V_{q , j} \cup \{ q \}$ for a unique pair
$(q , j)$ with $q \in I$ and $j \leq n(q)$.
\end{defn}
A moment's thought shows that if $\graph$ is a tree-like hypergraph
rooted at $p$ then the tree structure on $\supp (\graph)$ satisfying
the definition is unique (when $p$ is specified as the root).
Denote this tree by $\tree_p (\graph)$.

Sometimes it will be desirable to allow singleton hyperedges
(hyperedges consisting of a single vertex, $p$).
Rather than change the definitions, we introduce the notion
of a \Em{marked hypergraph}.  This is just a pair $(\graph , \marked)$,
where $\graph$ is a finite hypergraph and $\marked$ is any
subset of $\supp (\graph)$.  We think of $\marked$ as telling
us (by marking) which singleton edges $\{ p \}$ have been added to $\graph$.
Hypergraphs $\graph$ and $\graph'$ are defined to be isomorphic
if there is a bijection $\phi : \supp (\graph ) \to \supp (\graph')$
inducing a bijection at the level of hyperedges.  Marked
hypergraphs $(\graph , \marked)$ and $(\graph' , \marked')$
are isomorphic if $\phi$ can be chosen so that also $\phi (\marked)
= \marked'$.

In what follows, we will require a notion of weak convergence
of probability measures on hypergraphs and marked hypergraphs,
which in turn requires a metric on the space of marked hypergraphs
on the vertex set $\R$ rooted at $p$ (and we will re-normalize,
replacing prime $p$ by the real number $\rho=\rho_p:=p/y$, 
which will thus lie in the fixed
interval $(1,M]$).  It will turn out that
all but a vanishing fraction of our hypergraphs are tree-like,
so we need only to define the metric on tree-like hypergraphs
(e.g., by convention we take the distance between hypergraphs
to be $+\infty$ if either one is not tree-like).  If $\graph$
and $\hyp$ are two tree-like hypergraphs, define the distance
to be $+\infty$ if the two hypergraphs are not isomorphic,
and otherwise define the distance to be the least $\ee > 0$
such that there is a bijection $\phi : \supp (\graph) \to
\supp (\hyp)$ inducing an isomorphism on the hypergraphs,
and satisfying $|\phi (\rho) -\rho| \leq \ee$ for all $\rho \in \supp (\graph)$.
(Here we are dealing with re-normalized values of $p$, that is $\rho_p=p/y$, which are bounded.)
In other words, the topology is discrete on the graph structure
along with the product topology on the names of the vertices.
Formally,
$$d(\graph,\hyp) := \min_\phi \left \{ \max_{\rho \in \supp (\graph)}
   |\phi (\rho) - \rho|  \; : \; \phi \mbox{ is an isomorphism from }
   \supp (\graph) \mbox{ to } \supp (\hyp) \right \}  \, .$$
Define the distance between marked hypergraphs similarly, with
$\phi$ now restricted to isomorphisms of the marked hypergraphs.
Let $\mu$ and $\mu'$ be two probability measures on the space
of hypergraphs on the vertex set $\R$.  Say that a random pair
$(\graph , \graph')$ of hypergraphs is a coupling of $\mu$ and $\mu'$
when $\graph$ has law $\mu$ and $\graph'$ has law $\mu'$.  Define
the distance $d (\mu , \mu')$ between the probability measures $\mu$
and $\mu'$ to be the infimum of values $\ee > 0$ such that
there is a coupling $(\graph , \graph')$ of $\mu$ and $\mu'$
for which the probability of $d(\graph , \graph') > \ee$ is
at most $\ee$.  This is a standard metrization of the weak
topology, that is, $d (\mu_n , \mu) \to 0$ if and only if
$\int f \, d\mu_n \to \int f \, d\mu$ for all bounded and
weakly continuous functions $f$.

\subsection{The random hypergraph $\graph$ of $(M y)$-smooth numbers}
\label{ss:4.2}

Before we get started, here are a few words on notation.  As before,
we are selecting random positive integers $\leq x$, with $y (x)$ and $J_0 (x)$
as in Section~1.  Also, as before, we will choose an integer
$J := \lfloor \eta J_0 \rfloor$ for some $\eta > 0$.
We will choose a real $M > 1$ and keep track
of large prime factors in the interval $(y , M y)$.
By the term \Em{large prime}, we will mean a prime
in the interval $(y , M y)$.
We will also choose an integer $k \geq 1$ and keep track
only of numbers with at most $k$ large prime factors (factors
in the interval $(y , M y)$); we may even choose $k = \infty$ in the range implied 
by the limitations given to the uniformity of (\ref{eq:unif}).  We will
also specify an integer $m \geq 1$ which is interpreted as
the maximum chain length our algorithm will exploit when counting
pseudosmooths, where a chain is a sequence $a_1 , a_2 , \ldots , a_r$,
$r \leq m$, such that each consecutive pair $a_i , a_{i+1}$ share
a large prime factor $p_i \in (y , M y)$.  The first mission of this
subsection is to define a random hypergraph which will depend
on $M, J, x, m, k$ and a large prime $p \in (y , M y)$.
The full notation for this will be $\graph_{m,k,p}^{M,J,x}$.
However, in most of the results and constructions that follow,
$k, M$ and $J$ are fixed and $x$ is a size parameter fixed
during each construction, while $m$ and $p$ are dynamic
(the constructions are recursive in $m$ and $p$ and the
proofs inductive).  Because of this, we often reduce clutter
in the notation by writing simply $\graph_{m,p}$ with the
other four parameters understood.  In many of our lemmas,
arises the phrase, ``$f = o(1)$ as $x \to \infty$, uniformly
as $M$ and $\eta$ vary over bounded intervals and
$y < p < M y$.''  To be precise about this once and
for all, it means that there is a function $g$, going to
zero as $x$ goes to infinity, such that
$f (M , J , x, m , k , p) < g (M_0 , \eta , x , m , k)$
for all $M \leq M_0 , J \leq \eta J_0$ and $y < p < M y$
as $x \to \infty$.  This holds for any fixed $m, k, M_0 , \eta$.
Several times in Section~\ref{ss:conv} below we prove weak
convergence results.  Note: such  convergence results needing to be uniform, in the
manner just described, was the reason for metrizing the weak topology.

Now we move on to the constructions.
Fix an integer $x > 0$ and let $(\Omega_x , \F_x , \P_x)$ be a
probability space on which is defined a sequence $\{ X_1 , X_2 , \ldots \}$
of IID random variables whose common distribution is uniform
on the set $\{ 1 , 2 , \ldots , x \}$.  Let $y = y_0 (x)$ and
$J_0 (x) = x \pi (y) / \psi (x,y)$ be as in Section~1.
For each real $M > 1$ and each integer
$J > 0$, we will define a random hypergraph on the space
$(\Omega_x , \F_x , \P_x)$, which we will denote by $\graph^{M,J,x}$.

Given a real number $M > 1$, we keep track of prime factors
up to $M y$ as follows.  For any integer $X$ that
is $(M y)$-smooth, define the class $[X]$ to be the set of
primes $p$ for which $y < p < M y$ and $X$ is divisible by $p$
to an odd power, that is $p \in [X]$ if and only if $y < p < M y$ and
$p^i \, | \, X$ but $p^{i+1} \not| \; X$ for some odd integer $i$.
If $X$ is $y$-smooth, we define $[X]$ to be
the empty set.  If $X$ is not $(M y)$ smooth, we pick a symbol
(for probabilists, the traditional symbol is $\cemetery$) and set
$[X] = \cemetery$.

Now we define a random hypergraph with vertices in $\R^+$ by
$$\graph := \graph^{M,J,x} := \left \{ [X_j] \: : \; [X_j] \neq
   \cemetery \mbox{ and } \# [X_j] \geq 2
   \right \}_{1 \leq j \leq J} \, .$$
We remark that for a fixed $x$, the random hypergraphs $\graph^{M , J , x}$
are defined simultaneously for all $M$ and $J$.  In case it seems
strange to take $V = \R^+$ instead of $\Z^+$, it is because we will
be taking scaling limits.  Some easy but useful estimates are as follows.
\begin{proposition} \label{pr:counting}
Fix $M > 1$ and $\eta > 0$.  Let $J = \lfloor \eta J_0 \rfloor$
and let $[X_1] , [X_2] , \ldots$ and $\graph$ denote the random
variables on $(\Omega_x , \F_x , \P_x)$ constructed above.
For any finite set $S$ of primes, let
\[ N(S) = \{j : j \leq J; [X_j] = S\} \,.\]
\begin{enumerate}
\item For any finite set
$S$ of primes in $(y , M y)$ with $|S| \geq 2$, the number $N (S)$
 has asymptotic mean
\begin{equation} \label{eq:means}
\E_x N(S) \sim \eta \, \frac{y (\log y)^{|S|-1}}{\prod_{p \in S} p} \, .
\end{equation}
An upper bound, with an extra factor, is valid for all $S$:
\begin{equation} \label{eq:extra factor}
\E_x N(S) \leq 2^{|S|+1} \, \eta \,
   \frac{y (\log y)^{|S|-1}}{\prod_{p \in S} p} \, .
\end{equation}
\item For any set $\cT$ of hyperedges $S$, let $N(\cT) :=
\sum_{S \in \cT} N(S)$ denote the total number of hyperedges
in $\cT$.  Then, for any $\cT$, $\P_x (N(\cT) \geq 2)
\leq (\E_x N(\cT))^2$.
\item For any $p \in (y , M y)$,
the probability that there will be a prime $q \neq p$ such that
more than one hyperedge of $\graph$ 
contains both $p$ and $q$ goes to zero uniformly in $M \leq M_0$,
$\eta \leq \eta_0$ and $y < p , q \leq M y$.
\end{enumerate}
\end{proposition}

\noindent{\bf Proof.}
The means are computed by counting the number of
$a \leq x$ with $[a] = S$.  The number of integers of the form
$s \prod_{p \in S} p$ up to $x$ where $s$ is $y$-smooth is
$\psi (x / \prod_{p \in S} p , y)$.  The number of integers
of this form that are divisible by $q^2$ for some $q \in S$ is
bounded above by $\displaystyle{\sum_{q \in S} \psi 
\left ( \frac{x}{q \prod_{p \in S} p} , y \right )}$.  
This is easily shown to be asymptotically negligible compared to 
$B_S := \displaystyle{\psi \left ( \frac{x}{\prod_{p \in S} p} , y \right )}$ 
by (\ref{2.7}), using the fact that $\alpha$ remains bounded away from zero,
hence the number of $a \leq x$ with
$[a] = S$ is asymptotically equal to $B_S$.
By (\ref{eq:unif}), and using $\pi (y) \sim y / \log y$, we then have
\begin{eqnarray*}
\E_x N(S) & \sim & J \frac{\psi ( x / \prod_{p \in S} p , y)}{x} \\[1ex]
& \sim & \eta \frac{y (\log y)^{|S|-1}}{\prod_{p \in S} p}\,,
\end{eqnarray*}
which is (\ref{eq:means}).  Using (\ref{eq:unif unif}) instead
of (\ref{eq:unif}), and $\pi (y) \leq 2 y / \log y$ instead of
$\pi (y) \sim y / \log y$, gives (\ref{eq:extra factor}).

The second statement follows because $N(\cT)$ has a binomial
distribution.  For the third statement, let $H(p)$ denote the
event that there is some $q$ for which more than one hyperedge
arises containing $p$ and $q$.  Fix any primes $p_1 \neq p_2$.
Let $\cT_k$ denote the set of sets $S = \{ p_1 , p_2 , \ldots , p_k \}$
of distinct primes between $y$ and $My$ and let $\cT =
\cup_{k \geq 2} \cT_k$.  By the second statement of this proposition,
an upper bound for $H(p_1)$ may be obtained by summing any upper bound
for $(\E_x N(\cT))^2$ as $p_2$ ranges over primes between $y$ and $My$.
We compute this by bounding $\E_x N(\cT_k)$, then summing over $k$,
squaring, and summing over $p_2$.  Thus we begin by
using (\ref{eq:extra factor}) with $S = \{ p_1 , p_2 , \ldots , p_k \}$
to obtain
$$\E N(S) \leq 2^{k+1} \eta y (\log y)^{k-1} \prod_{p \in S}
   \frac{1}{p} \, .$$
Summing this over all choices of $p_3 , \ldots , p_k$
and using (\ref{eq:PNT 2}) for the last inequality then gives
\begin{eqnarray*}
\E N(\cT_k) & \leq & \frac{2^{k+1} \eta \, y \log y}{p_1 p_2}
   \sum_{p_3 < \cdots < p_k} \prod_{j=3}^k \frac{\log y}{p_j} \\[1ex]
& \leq & \frac{2^{k+1} \eta \, y \log y}{p_1 p_2} \frac{1}{(k-2)!}
   \sum_{p_3 , \cdots , p_k} \prod_{j=3}^k \frac{\log y}{p_j} \\[1ex]
& \leq & \frac{2^{k+1} \eta \, y \log y}{p_1 p_2} \frac{1}{(k-2)!}
   \prod_{j=3}^k (2 \log M) \, .
\end{eqnarray*}
We sum this over all integers $k\geq 3$ so that
$$\E N(\cT) \leq \frac{8 M^4 \eta y \log y}{p_1 p_2} \leq
\frac{8 M_0^4 \eta_0   \log y}{  p_2}\, ,$$
since $y/p_1<1$.  Squaring, noting that $1/p_2<1/y$ and $\log y<\log p_2$, 
we obtain a quantity bounded above by a constant multiple of
$$\frac{\log y}{y} \sum_{y < p_2 \leq My} \frac{\log y}{p_2}$$
By~(\ref{eq:PNT 1}) this is $\displaystyle{O(\frac{\log y}{y})}$;
this completes the proof, as we only needed to show $o(1)$.  
$\Cox$

We now define sub-hypergraphs $\graph_{m,p}$ of the random hypergraph
$\graph$, culled so as to be tree-like and rooted at $p$.  They
are deterministic functions of the variables $X_1 , \ldots ,
X_{J}$, and they will bear witness to the creation of
pseudo-smooth numbers.  They depend on the parameters $M , J, x$
and $k$, which are fixed throughout the construction and suppressed in
the notation.  We remark that the definition makes sense for $k = \infty$.
\begin{defn}[The sub-hypergraph $\graph_{m,p}$ and marked set
$\marked_{m,p}$] \label{def:G}
We define hypergraphs $\graph_{m,p} (j)$ recursively for $m \geq 1$
and $1 \leq j \leq J$ as follows.
\begin{itemize}
\item Let $T_0 (p) := \{ p \}$ and $\graph_{0,p} := \emptyset$, 
taking $\supp (\graph_{0,p}) = \{ p \}$ by convention.
\item For each $m \geq 1$, define $\graph_{m,p} (0) := \graph_{m-1,p}$.
For $j \geq 1$, define $\graph_{m,p} (j) :=
\graph_{m,p} (j-1) \cup \{ [X_j] \}$ if $[X_j]$ intersects
$\supp (\graph_{m,p} (j-1))$ in a single element of
$T_{m-1} (p)$ and $2 \leq |[X_j]| \leq k$.  Otherwise,
let  $\graph_{m,p} (j) := \graph_{m,p} (j-1)$.
Define $\graph_{m,p} := \graph_{m,p} (J)$.  Define
$T_m (p) := \supp (\graph_{m,p}) \setminus \supp (\graph_{m-1,p})$.
\end{itemize}
Let $\marked$ denote the set of primes $q$ with $y < q < M y$
such that $[X_j] = q$ for some $j \leq J$.  Let $\marked_{m,p} :=
\marked \cap \supp (\graph_{m,p})$.  Then $(\graph_{m,p} ,
\marked_{m,p})$ is a marked sub-hypergraph, which we will use later
to witness the creation of pseudo-smooths.
\end{defn}
%
Informally, $\graph_{1,p}$ takes all hyperedges of $\graph$
that contain $p$ except for those creating a collision (that is, 
a cycle on hyperedges), using the order in which they were 
generated to settle collisions.
Then, $\graph_{2,p}$ starts over, taking all hyperedges containing
each of the vertices added in the previous step, except for those
that cause collisions.  In the end, the list of hyperedges is
swept through, in order, $m$ times.  The informal interpretation of
$T_m (p)$ is the set of primes that first appear at distance $m$
from $p$ in our tree-like hypergraph; the informal interpretation of
$U_{m,p}$ is the set of primes within distance $m$ of $p$ that
appear as hyperedges of cardinality one.  
\begin{lemma} \label{lem:moments}
For any $\eta, M, x$ and $p$,
$$\E_x \left | \graph_{1,p} \right | \leq (2M - 1) \eta \frac{y}{p} \, .$$
\end{lemma}

\noindent{\bf Proof.} By construction, the hypergraph
$\graph_{1,p}$ is a subset of the restriction of $\graph$
to hyperedges containing $p$.  Therefore,
$$\E_x |\graph_{1,p}| \leq \sum_S \E_x N(S)$$
where the sum is over such sets $S$.  Break down the sum
by the cardinality of $S$.  The sum over $|S| = k$ is
$1 / (k-1)!$ times the sum over ordered sets of primes $p = p_1 ,
p_2 , \ldots , p_k$ in the range $(y , M y)$.  The sum over
ordered such sets is bounded above by the sum over ordered
$k$-tuples in which repetition is allowed.  Thus
$$\E_x |\graph_{1,p}| \leq \sum_{k \geq 2} \frac{1}{(k-1)!}
   \sum_{p_2 , \ldots , p_k} \E_x N(p , p_2 , \ldots , p_k)$$
where the summand is zero, by convention, if there is a repetition.
When there is no repetition, we obtain an estimate from (\ref{eq:means}),
which implies the upper bound
$$\E_x |\graph_{1,p}| \leq \sum_{k \geq 2} \frac{1}{(k-1)!}
   \sum_{p_2 , \ldots , p_k} \eta \frac{y}{p} \prod_{i=2}^k
   \frac{\log y}{p_i} \, .$$
The inner sum factors as a power, yielding
$$\E_x |\graph_{1,p}| \leq \eta \frac{y}{p} \,
   \sum_{k \geq 2} \frac{1}{(k-1)!}
   \left ( \sum_{y < q < M y} \frac{\log y}{q}
   \right )^{k-1} \, .$$
By the prime number theorem, $\sum_{y < q < M y} (\log y)/q \to \log M$,
and is never more than $\log (2M)$, whence
$$\E_x |\graph_{1,p}| \leq \eta \frac{y}{p} \,
   \sum_{k \geq 2} \frac{(\log (2M))^{k-1}}{(k-1)!}
   = (2M - 1) \eta \frac{y}{p} \, . $$
$\Cox$

\begin{corollary} \label{cor:size}
$$\E_x |\graph_{1,p}|^2 \leq \E_x |\graph_{1,p}|
   + \left ( \E_x |\graph_{1,p}| \right )^2$$
and
$$\E_x |\graph_{m,p}| \leq \left (1 + 2 \eta M \right )^m \frac{y}{p} \, .$$
\end{corollary}

\noindent{\bf Proof.} For the first statement, note that for $S \neq T$,
the events $\{ S \in \graph_{1,p} \}$ and $\{ T \in \graph_{1,p} \}$
are negatively correlated. (Recall that two events are negatively correlated, if the probability of their conjunction
is {\em at most} the product of the probabilities of the events.)  This is because the events $\{ [X_i] = S \}$
and $\{ [X_j] = T \}$ are independent,  unless $i=j$,  in which case
they are negatively correlated.  It follows that
\begin{eqnarray*}
\E_x |\graph_{1,p}|^2 & \leq & \sum_{S,T}
   \P_x \left ( S , T \in \graph_{1,p} \right ) \\
& \leq & \left ( \E_x |\graph_{1,p}| \right )^2 +
   \E_x |\graph_{1,p}| \, .
\end{eqnarray*}

\noindent

For the second statment, induct on $m$.  Conditional on $\graph_{m-1,p}$,
the random hypergraph $\graph_{m,p}$ is stochastically dominated by
the union of $\graph_{m-1,p}$ with a collection of hyperedges whose
conditional distribution given $\graph_{m-1,p}$ is described as follows:
for each $q \in T_{m-1} (p)$, and for each finite subset $S$ of
primes in $(y , My)$ containing $q$, the hyperedge $S$ is added
independently with probaiblity $N(S)$.  By induction, the
mean number of such $q$ is at most $(1 + 2 \eta M)^{m-1} y / p$.
Bounding the mean of each Poisson  variable from above by $2 \eta M$,
we complete the induction.   $\Cox$

\

The number of pseudo-smooths generated by time $j$, by definition,
is the difference between $j$ and the $\FF_2$-rank of the collection
$[X_1] , \ldots , [X_j]$, made into a $\FF_2$-vector space
by using the symmetric difference operation $[X_i] \symdif [X_j]$.
To count this, we count the number of $j$ for which $[X_j]$ is in
the $\symdif$-span of $[X_1] , \ldots , [X_{j-1}]$, which we denote  
by $\langle [X_1] , \ldots , [X_{j-1}] \rangle$.  This includes
the case where $[X_j] = \emptyset$ ($y$-smooth numbers), $[X_j]
= [X_i]$ for some $i < j$ (the one large prime case), as well as
more complicated cases.  It turns out that not much is lost if
we include only one more class of cases.  For each prime $p$ in
the interval $(y , M y)$, and each positive integer $j$,
we define an event $\chi_{m,k,p}^{M,j}$ whose informal interpretation
is that $\{ p \}$ is in the span of $\langle [X_1 , \ldots , [X_j] \}$
and that this fact is witnessed by classes $[X_i]$ of cardinality
at most $k$, having indices $i \leq j$.  A proposition immediately
following the definition verifies the interpretation.  The parameters
$k,x,j$ and $M$ will now be fixed throughout the definition and
suppressed from the notation.
\begin{defn}[$\chi$ for general marked rooted trees]\label{def:chi-marked}
~\\ \vspace{-0.25in}
\begin{enumerate}
\item Let $(G , U)$ be any marked hypergraph rooted at
a vertex $p$.  For $q \in \supp (G)$, define the height $\ell (q)$
to be the length of the longest non-backtracking path from $q$ to
the leaves of $G$, or more accurately, of the tree $\tree_p (G)$.
\item Define an event $\chi (q) = \chi (G,U,q)$ by recursion on $\ell (q)$.
If $\ell (q) = 0$, define the event $\chi (q)$
to hold if and only $q \in U$.  If $\ell (q) > 0$,
let $r$ denote the distance from $p$ to $q$ in $\tree_p (G)$ and
define $\chi (q)$ to hold if and only if there is some hyperedge $S \in G$
such that~$(i)$~$S \subseteq T_{r+1} (p) \cup \{ q \}$ (that is, $S$ is
a hyperedge that appears first at distance $r+1$ from $p$, and is a ``child'' of $q$), and~$(ii)$
the event $\chi (q')$ occurs for each $q' \in S$ other than $q$.
\item Finally, let $\chi (G,U)$ denote the event $\chi (G,U) (p)$.
\end{enumerate}
\end{defn}
\begin{unremarks}
Note that the recursion is well founded because $\ell (q')
\leq \ell (q) - 1$ for all such $q'$.  Also note that
in the recursive part of the definition, we allow
$S$ to equal $\{ q \}$, in which case~$(ii)$ is vacuously
satisfied.
\end{unremarks}
\begin{defn}[smooth primes witnessed in an $m$-neighborhood] \label{def:chi}

\

\noindent
If $\graph_{m,p}$ is not tree-like, we define $\chi_{m,p}$
not to occur.  If $\graph_{m,p}$ is tree-like, we define
$\chi_{m,p} (q) := \chi_{m,p} (\graph_{m,p} , \marked_{m,p} , q)$,
whence,
$$\chi_{m,p} := \chi (\graph_{m,p} , \marked_{m,p}, p) \, .$$
\end{defn}
Let $\VS$ denote the vector space over $\FF_2$ whose basis is the
set of symbols
$$\{ \delta_p : p \mbox{ is a prime and } y < p < My \}.$$
Identify each class $[X]$ with the element $\sum_{p \in X} \delta_p$
of $\VS$.  In the following proposition, $\langle [X_1] , \ldots , [X_j]
\rangle$ denotes the span of $\{ [X_1] , \ldots , [X_j] \}$ in $\VS$.
\begin{proposition} \label{pr:branching}
For any $m \geq 1$, the event $\chi_{m,p} (q)$ implies $\{ q \} \in
\langle [X_1] , \ldots , [X_j] \rangle$.  In particular,
$$\chi_{m,p} \;  \Longrightarrow  \;
   \{ p \} \in \langle [X_1] , \ldots , [X_j] \rangle \, .$$
\end{proposition}

\noindent{\bf Proof.} By induction on $\ell (q)\geq 0$.  If $\ell (q) = 0$
then $\chi_{m,p} (q)$ implies $[X_j] = \{ q \}$ for some $j \leq J$,
which immediately implies $\{ q \} \in \langle [X_1] , \ldots , [X_j]
\rangle$.  Now suppose $\ell (q) \geq 1$.  If $\chi_{m,p} (q)$
holds, let $j$ satisfy $(i)$ and $(ii)$ of the definition with $q=p$.
For each $q' \in [X_j]$ distinct from $q$, $\ell (q') \leq
\ell (q) - 1$, whence by induction, $\{ q' \} \in
\langle [X_1] , \ldots , [X_j] \rangle$ for all such $q'$.
This, along with the trivial observation that $[X_j] \in
\langle [X_1] , \ldots , [X_j] \rangle$, implies $\{ q \}
\in \langle [X_1] , \ldots , [X_j] \rangle$, which completes
the induction.
$\Cox$

\

It follows from this that for any $m$, the number of linear
dependences among \\
$\{ [X_1] , \ldots , [X_J] \}$ is bounded
from below by
\begin{equation} \label{eq:count}
\# \{ j \leq J : \mbox{ for all } p \in [X_j], \mbox{ the singleton } 
   \{ p \} \mbox{ is in the span } \langle [X_1] , \ldots , [X_{j-1}] 
   \rangle \} \, .
\end{equation}

\subsection{Construction of the limit object $\hyp_m$}
\label{ss:4.3}

An informal description of the limit object is as follows.
The root, $\rho$, gets hyperedges $\{ \rho , \rho_1 , \ldots , \rho_k \}$
independently, with the probability of such a hyperedge arising
in a small volume element $\{ \rho \} \times [\rho_1 , \rho_1 + d\rho_1]
\times \cdots \times [\rho_k , \rho_k + d\rho_k]$ equal to 
$$\frac{d\rho_1 \cdots d\rho_k}{\rho \rho_1 \cdot \rho_k} \, .$$
Recursively, for $m$ iterations, each vertex newly added in the last 
iteration gets new hyperedges in the same way.

Formally, the limit object is best described in terms of Poisson 
processes.  We briefly summarize definitions and properties of these,
referring the reader to~\cite{poisson} for further details.
Given a measure space $(\SS , \B)$ with a $\sigma$-finite
measure $\mu$, a Poisson process with intensity $\mu$ is
a collection of random variables $\{ N(S) = N(S)(\omega) : S \in \B \}$
on some probability space $(\Omega , \F , \P)$ satisfying
the following properties:

(1) Countable additivity in $S$:
if $\A$ is a collection of disjoint elements of $\B$ then
$N \left( \bigcup_{S \in \A} S\right) = \sum_{S \in \A} N(S)$;

(2) Counting measure: $N(S)$ takes values in the nonnegative integers;

(3) Poisson distribution: for fixed $S$, the random variable $N(S)$
is distributed as a Poisson distribution with mean $\mu (S)$;

(4) Independence:
if $S,T$ are disjoint elements of $\B$ then $N(S)$ and $N(T)$ are
independent.

\noindent A number of constructions are available to prove the
existence of such a process.  

If $\mu$ is nonatomic, then with
probability~1, the random counting measure $N$ gives measure
at most~1 to every point $s \in \SS$.  It follows that the
random measure $N(S)$ is the sum of point masses $\delta_s$,
as $s$ ranges over some finite or countable subset of $\SS$;
we denote this set by $\supp (N)$ and refer to $\supp (N)$
as ``the points of the Poisson process''.  The cardinality
of $\supp (N)$ is a Poisson random variable with mean $\mu (\SS)$.

Fix a real number $M > 1$.  Fix also a real $\eta > 0$ and an integer
$k \geq 2$.  We construct a random hypergraph $\hyp_{m,p} =
\hyp_{m,k,p}^{M,\eta}$ on a new probability space $(\Omega , \F , \P)$
whose vertex set is the real interval $[1,M]$.  The collection
$[1,M]_j$ of subsets of $[1,M]$ of cardinality $j$ may be
identified with the sector $W_j \subseteq \R^j$ defined by
$$W_j := \{ (\rho_1 , \ldots , \rho_j) \in \R^j : 1 \leq \rho_1 < \cdots
   < \rho_j \leq M \} \, .$$
Let $d\pp / (\rho_1 , \ldots , \rho_j)$ denote the image under
this identification of the measure whose density with respect
to Lebesgue measure is $1 / (\rho_1 \cdots \rho_j)$.  Observe that
the total mass of the measure $d\pp / (\rho_1 \cdots \rho_j)$ is given
$(\log M)^j / j!$.  Now define a measure $\mu_k$ on the union
$\bigcup_{j=1}^k [1,M]_j$ by $\mu_k = \sum_{j=1}^k
d\pp / (\rho_1 \cdots \rho_j)$.  Let $\mu$ denote the increasing limit
of $\mu_k$ as $k \to \infty$.  We see that $\mu$ has finite total mass:
$$||\mu|| \, = \, \sum_{j=2}^\infty \frac{(\log M)^j}{j!} \, = \,
   M - 1 - \log M \, .$$
Fix $\rho \in [1,M]$ and define an operation $\sigma_\rho$ by
$\sigma_\rho (S) = S \cup \{ \rho \}$.  Define the measure $\mu_k^{+\rho}$ by
$\mu_k^{+\rho} = \mu_k \circ \sigma_\rho^{-1}$.  In other words, $\mu_k^{+\rho}$
is the measure corresponding to ``choosing a set according to $\mu_k$''
and then adding the element $\rho$.  (Here the quotes are  to remind the reader that the finite measure $\mu_k$ is not a probability measure).  Thus all the measures
$\mu_k^{+\rho}$ as well as the increasing limit $\mu^{+\rho}$ are
supported on finite sets of cardinality at least~2.

Let $\tau \in [1,M]$ (here $\tau$ plays the role of $q/y$, just as
$\rho$ plays the role of $p/y$).  Let $\nu_\tau = \nu_{k,\tau}^{M ,\eta}$ (as usual,
we suppress quantities that are, for the moment, fixed) be the law
of the points of a Poisson process with intensity $\mu^{+\rho} / \tau$.
Observe that each point of the process is a finite subset $S$
of $[1,M]$ with $\rho \in S$.  Because the intensity measure has
finite mass, the law of the set of points is the law of a random
finite set of hyperedges $S \subseteq [1,M]$.  By non-atomicity
of Lebesgue measure, we see that with probability~1, this is a
tree-like hypergraph rooted at $\rho$, all of whose hyperedges
contain $\rho$.

\begin{defn}[The marked graph $(\hyp_{m,\rho} , \vmarked_{m,\rho})$]
\label{def:hyp}
We now construct the random hypergraphs $\hyp_{m,\rho} =
\hyp_{m,k,\rho}^{M , \eta}$, by recursion on $m$.  For $m=1$,
choose $\hyp_{1,\rho}$ from the law $\nu_\rho$.
For $m \geq 1$, let $T_{m,\rho} = \supp (\hyp_{m,\rho})
\setminus \supp (\hyp_{m-1,\rho})$, taking $\supp (\hyp_{0,\rho})
= \{ \rho \}$ by convention.  For the recursion step,
choose random hypergraphs $\hyp_{m , \tau}$ independently from
respective laws $\nu_\tau$, as $\tau$ varies over $T_{m,\rho}$, and let
$\hyp_{m+1,\rho}$ be the union of $\hyp_{m,\rho}$ with all
the sets $\hyp_{m+1 , \tau} \,$.  It is again immediate that
each $\hyp_{m,\rho}$ is tree-like.  Finally, we define a set
of marks $\vmarked_{m,\rho}$\,, by choosing each $\tau \in \supp (\hyp_{m,\rho})$
independently, with probability $1 - e^{- \eta / \tau}$.
\end{defn}

Now, using Definition~\ref{def:chi} once more, define events
\begin{eqnarray*}
\chi_{m,\rho}' (\tau) & := & \chi (\hyp_{m,\rho} , \vmarked_{m,\rho} , \tau) \, ;  \\
\chi_{m,\rho}' & := & \chi (\hyp_{m,\rho} , \vmarked_{m,\rho})  \, .
\end{eqnarray*}
These are events on the space $\Omega$ analogous to the events
$\chi_{m,p} (q)$ and $\chi_{m,p}$ defined on the space $\Omega_x$.
Denote
$$\theta_m (\rho) := \theta_{m,k}^{M,\eta} (\rho) := \P (\chi_{m,\rho}') \, .$$

\subsection{Convergence of $\graph$ to $\hyp$, and consequently,
of $\P_x (\chi)$ to $\theta$}  \label{ss:conv}
\label{ss:4.4}

In this subsection we prove convergence results which will be used
to compute the rate of accumulation of pseudo-smooth numbers.
\begin{theorem} \label{th:chi}
Fix integers $m,k \geq 1$ and any real $M > 1$.  Then
\begin{equation} \label{eq:chi}
\P_x (\chi_{m,k,p}^{M,j}) = (1 + o(1)) \theta_{m,k}^{M,j/J_0} (p/y)
\end{equation}
uniformly as $p$ varies over primes in the interval $(y , M y)$
and $j/J_0$ remains bounded.  More generally, for any $r \geq 1$
and any $p_1 , \ldots , p_r$,
\begin{equation} \label{eq:chi square}
\P_x \left ( \bigcap_{i=1}^r \chi_{m,k,p_i}^{M,j} \right ) = (1 + o(1))
   \prod_{i=1}^r \theta_{m,k}^{M,j / J_0} (p_i/y)\,,
\end{equation}
uniformly as $p_1 , \ldots , p_r$ vary over primes in the
interval $(y , M y)$.
\end{theorem}

The proof of this theorem is essentially to show that the
rescaled random graph $y^{-1} \graph_{m,p}$ converges
weakly to $\hyp_{m , p/y}$.  We encapsulate what we need
in the following lemmas.  All of these are routine Poisson
convergence lemmas.  In each case, the lemmas hold for any
fixed $k$, and with  $k = \infty$ in the range of uniformity
given for  (\ref{eq:unif}).   
\begin{lemma} \label{lem:conv 1}
As $x \to \infty$, the distance in the weak metric between the
random hypergraph $\graph_{1,p}^{M,j,x}$ and the random hypergraph
$\hyp_{1,p/y}^{M,j/J_0}$ goes to zero, uniformly as $M$ and
$j/J_0$ vary over bounded intervals and $y < p < M y$.
\end{lemma}

\noindent{\bf Proof.}
As a preliminary computation, let $\graph_{1,p}'$ denote the
subset of $\graph$ of all hyperedges containing $\{ p \}$.
We claim that $\P (\graph_{1,p} = \graph_{1,p}') \to 1$.
Indeed, the complementary event requires that a collision
occur, entailing two hyperedges both to contain $\{ p \}$
and $\{ q \}$ for some $q$.  By the last part of
Proposition~\ref{pr:counting}, this probability goes to
zero uniformly (and even for $k = \infty$
in the range allowed by using (\ref{eq:unif})).

Next, let $\Xi = (\tau_1, \tau_1'] \times \cdots \times (\tau_n, \tau_n']$ be any rectangular
subset of the sector $\sector_n$ and let $\Xi_x$ denote the set of
sets, $S$, of $n$ primes, each between $y$ and $My$,
such that $y^{-1} S \in \Xi$.  As in Proposition~\ref{pr:counting}, let $N(\sigma_p (\Xi_x))$ denote
the number of $j \leq J$ such that $[X_j] \in \sigma_p (\Xi_x)$.
Using (\ref{eq:means}), we estimate
\begin{eqnarray*}
\E_x N(\sigma_p(\Xi_x)) & = & \sum_{S \in \sigma_p (\Xi_x)} \E_x N(S)
   \\[1ex]
& \sim & \sum_{S \in \sigma_p (\Xi_x)} \eta \frac{y (\log y)^n}
   {p \, \prod_{q \in S} q} \, .
\end{eqnarray*}
Factoring the sum of products gives the equivalent expression
$$\E_x N(\sigma_p (\Xi_x)) \sim \eta \frac{y}{p} \prod_{i=1}^n
   \sum_{\tau_i y < q \leq  \tau_i' y} \frac{\log y}{q} \, .$$
By the prime number theorem, this converges to $\nu_{p/y} (\Xi)$.

Finally, let us see that $y^{-1} \graph_{1,p}$ converges to a
Poisson process with intensity $\nu_\rho$ where $\rho=p/y$; by construction, this
is the distribution of $\hyp_{1,\rho}$, and therefore this will
complete the proof of the lemma.  We need to show that for
any disjoint sets $\Xi^{(1)} , \ldots , \Xi^{(n)}$, the
respective numbers $N^{(i)}$ of hyperedges in $y^{-1} \graph_{1,p}$
in $\Xi^{(i)}$ converge in disribution to independent Poissons
with means $\nu_\rho (\Xi_i)$.  It suffices to prove this for
$\graph_{1,p}'$ in place of $\graph_{1,p}$ because we have
seen these are equal with probability $1 - o(1)$.

We have already verified that the means are $\nu_\rho (\Xi^{(i)})$
when $\Xi^{(i)}$ are rectangles, which implies the same
result for all measurable $\Xi$.  To obtain the
joint Poisson distribution, it is easiest to Poissonize.  Replace
$\graph_{1,p}'$ by $\graph_{1,p}''$, defined identically to
$\graph_{1,p}'$ except with $J$ replaced by a Poisson variable
$J'$ of mean $J$.  For this random graph, the numbers $(N^{(i)})''$
of hyperedges of $\graph_{1,p}''$ in the rescaled $\Xi^{(i)}$
are exactly independent Poissons with the given means.  The
key observation is that
$$\P_x (\graph_{1,p}' \neq \graph_{1,p}'') = O(J_0^{-1/2}) \, .$$
To see this, note that $\E_x |J' - J| = O(\sqrt{J_0})$.  Therefore,
\begin{equation} \label{eq:poissonize}
\P_x (\graph_{1,p}' \neq \graph_{1,p}'')
   = O \left ( \sqrt{J_0} \; \P_x (p \in [X_1]) \right )
   = O \left ( J_0^{-1/2} \; \E_x |\graph_{1,p}'| \right )
   = O \left ( J_0^{-1/2} \right )\,,
\end{equation}
by Corollary~\ref{cor:size}.
$\Cox$

\begin{lemma} \label{lem:conv 2}
As $x \to \infty$, the distance in the weak metric between the
$n$-tuple of random hypergraphs
$$y^{-1} \left ( \graph_{1,p_i}^{M,j,x} \right )_{1 \leq i \leq n}$$
and the product of the laws of the hypergraphs
$\hyp_{1,p_i/y}^{M,j/J_0}$ goes to zero, uniformly as $M$ and
$j/J_0$ vary over bounded intervals and $y < p_i < M y$.
\end{lemma}

\noindent{\bf Proof.}
This is the same proof with only one difference, as follows.
To check that $\graph_{1,p_i} = \graph_{1,p_i}'$ with probability
tending to~1, one observes that~(3) of Proposition~\ref{pr:counting}
holds simultaneously for $p_1 , \ldots , p_n$.  All else is the same,
once one observes that Poissonization gives (\ref{eq:poissonize})
simultaneously for all $p_1 , \ldots , p_n$.
$\Cox$

\begin{lemma} \label{lem:conv 3}
As $x \to \infty$, the distance in the weak metric between the
random hypergraph $y^{-1} \graph_{m,p}^{M,j,x}$ and the random hypergraph
$\hyp_{m,p/y}^{M,j/J_0}$ goes to zero, uniformly as $M$ and
$j/J_0$ vary over bounded intervals and $y < p < M y$.
Similarly, the distance between the law of the random $n$-tuple
$\disp{y^{-1} (\graph_{m,p_i}^{M,j,x})_{1 \leq i \leq n}}$
and the product of the laws of $\hyp_{m,p_i/y}^{M,j/J_0}$
goes to zero with the same uniformity in $M,j/J_0$ and $\{ p_i \}$.
\end{lemma}

\noindent{\bf Proof.}
We induct on $m$.  For $m=1$ this was shown in Lemma~\ref{lem:conv 1}.
Now let $m \geq 2$ and assume for induction that the result holds for
$m-1$.  If $\graph_{m,p}$ is tree-like, let $r := |T_1 (p)|$ and let
$G_1 , \ldots , G_r$ denote the subtrees of $\tree_p (\graph_{m,p})$
from the vertices $q_1 , \ldots , q_r$ of $T_1 (p)$.  Let
$\graph (1) , \ldots , \graph (r)$ denote the corresponding
hypergraphs, that is, $\graph (i)$ is the hypergraph rooted at
$q_i$ whose hyperedges are those of $\graph_{m,p}$ whose support
is a subset of the vertices of $G_i$.  We will show that the joint
conditional distribution of $y^{-1} (\graph (1) , \ldots , \graph (r))$
given $\graph_{1,p}$ converges to the product of the laws
of $\hyp_{m-1 , q_i/y}$.
By the recursive construction of $\hyp_{m,p/y}$ and the fact that
$\graph_{m,p}$ is tree-like with probability approaching~1, this
will complete the proof of the lemma.

Consider the hypergraph $\graph_{m-1 , q_i}'$.
If this is tree-like, let $H_i$ be the subtree obtained by
removing the unique hyperdege containing $p$ and $q_i$,
and restricting to the connected component rooted at $q_i$.
If these are disjoint for $1 \leq i \leq r$, then
$\graph (i) = H_i$ for each $i$.  The probability that
all the hypergraphs $\graph_{m-1,q_i}'$ are tree-like is
asymptotically~1.  The probability of a collision is bounded
above by
$$\sum_{y < q < M y} \sum_{1,j=1}^r \P_x \left ( q \in
   \supp (\graph_{m-1,q_i}' \cap \supp (\graph_{m-1,q_j}' \right ) \, .$$
The probability that $q \in \supp (\graph_{m-1,q_j})$, conditional
on $|\graph_{m-1,q_j}'|$, is $O(|\graph_{m-1,q_j}'| / \pi (y))$.
This is true as well for $q_i$, and the two events are independent.
Therefore, the probability of a collision is
$$O \left ( (\E_x r^2)   \frac{(\E_x |\graph_{m-1,y}|)^2}{\pi (y)}
   \right ) \, .$$
By Corollary~\ref{cor:size}, we obtain the upper bound $O(1 / \pi (y))$.

Next, we claim that the conditional distribution of $H_i$ given
$\graph_{1,p}'$ is asymptotically equal to the unconditional
distribution of $\graph_{m-1 , q_i}'$.  Indeed, $\graph_{1,p}'$
is measurable with respect to the $\sigma$-field generated by
the events $\{ S \in \graph : p \in S \}$.  This is independent
of the events $\{ S \in \graph : p \notin S \}$, so conditional
on $\graph_{1,p}'$, $H_i$ has the distribution of $\graph_{m-1,q_i}''$
where the double prime means that all hyperedges containing
$p$ were excluded at every step of the construction.  We already
know that $\graph_{m-1,q-1}''$ is asymptotically distributed as
$\graph_{m-1,q-1}'$, verifying the claim.  Moreover, the same
argument shows that the joint conditional law of $H_1 , \ldots , H_r)$
given $\graph_{1,p}'$ is asymptotically the product of the laws
for each $i \leq r$.

Finally, by the induction hypothesis, the unconditional distribution
of $\graph_{m_i , q_i}'$ is asymptotically that of $\hyp_{m-1 , q_i/y}$.
Therefore, since with probability approaching~1 all the graphs
$\graph_{m-1,q_i}'$ are tree-like and there are no collisions, we
have shown what we need.
$\Cox$

\begin{lemma} \label{lem:conv 4}
As $x \to \infty$, the distance in the weak metric between the
random marked hypergraph $y^{-1} (\graph_{m,p}^{M,j,x} , \marked_{m,p})$
and the random marked hypergraph
$(\hyp_{m,p/y}^{M,j/J_0} , \vmarked_{m,p/y})$ goes to zero,
uniformly as $M$ and $j/J_0$ vary over bounded intervals and
$y < p < M y$.  More generally, the distance between an
$n$-tuple of marked graphs
$$y^{-1} \left ( \graph_{m,p_i}^{M,j,x} , \marked_{m,p_i}
   \right )_{1 \leq i \leq n}$$
and the product of the laws of the random marked hypergraphs
$(\hyp_{m,p_i/y}^{M,j/J_0} , \vmarked_{m,p_i/y})$ goes to zero
uniformly as $M$ and $j/J_0$ vary over bounded intervals and
$y < p_1 , \ldots , p_n < M y$.
\end{lemma}

\noindent{\bf Proof.}
Observe that the conditional probabilities of
$q \in \marked_{m,p}$ given $\graph_{m,p}$ are independent
and given by $1 - e^{-\eta y / q}$ as $q$ varies over $\supp (\graph_{m,p})$.
This  is true since, in the limit ($x,y \to \infty$ and $J = \eta  x \pi(y)/\psi(x,y))$,
the events $|\{j : [X_j] = \{q_i\}, j=1, \ldots , J\}|$ for fixed $q_1, q_2, \ldots , q_r$ are independent
Poisson random variables with mean $\sim \eta y/q_i$.
And once it is known, in the limit, that the events $\{q\in U_{m,p}\}$ given
$\graph_{m,p}$, with $q$ running over $\supp(\graph_{m,p})$ are independent with probability 
$1 - e^{-\eta y/q}$, then the first part of the lemma is proven; the second part is analogous.

$\Cox$


\noindent{\bf Proof of Theorem~\ref{th:chi}.}
Begin with (\ref{eq:chi}).  For any marked graph $(G,U)$,
$\chi (G,U)$ depends only on the marked hypergraph structure
of $(G,U)$ and not the names of the vertices.  Because the
topology on graph structure is discrete, $\chi$ is continuous.
The weak topology on measure is characterized by convergence of
integrals of bounded continuous functions, so (\ref{eq:chi})
follows from the first conclusion of Lemma~\ref{lem:conv 4}.
For any fixed bounded continuous function, such as $\chi$,
the difference in the integrals is bounded as a function of
the distance bewteen the measures, whence the uniform convergence
in Lemma~\ref{lem:conv 4} transfers to the required uniform
convergence in (\ref{eq:chi}).  The proof of (\ref{eq:chi square})
is identical, using the $n$-tuple convergence in Lemma~\ref{lem:conv 4}
in place of convergence of the single marked hypergraph.
$\Cox$

\subsection{Computation of $\theta$}
\label{ss:4.5}

We begin by computing $\theta_m (\rho)$.  Recall the definition
of the functions $\gamma_{m,M,k} (u)$ in (\ref{eq:gamma}).
\begin{lemma} \label{lem:theta}
$$\theta_{m,k}^{M , \eta} (\rho) = 1 - e^{- \gamma_{m,M,k} (\eta) / \rho} \, .$$
\end{lemma}

\noindent{\bf Proof.} The quantities $M, \eta$ and $k$ will be
fixed throughout the proof, so we write $\theta_m$ for
$\theta_{m,k}^{M , \eta}$. The proof is by induction on $m$.
By definition, $1 - \theta_0 (\rho)$ is the probability that
$\rho \notin \vmarked_{m,\rho}$,  which is $e^{-\eta/\rho}$ by construction.
This establishes the result for $m=0$.

Now suppose that $m\geq 1$.  The set of hyperedges $S \in \hyp_{1,\rho}$
is, by construction, a Poisson process with intensity
$\nu_{\rho}$.  The complement of $\chi_{m,\rho}$ is
the intersection of $\rho \notin \vmarked_{m,\rho}$ with
the event that for all hyperedges $S \in \hyp_{1,\rho}$
of cardinality between~2 and~$k$,  there is some 
$\tau \in  S\setminus \{ \rho \}$  that is not in $ \vmarked_{m,\rho}$.  We have, by induction,
\begin{eqnarray}
1 - \theta_{m+1} (\rho) & = & e^{-\eta/\rho}\ \E \left [ \prod_{S \in \hyp_{1,\rho}}
   \left ( 1 - \prod_{\tau \in S \setminus \{ \rho \}} \theta_{m} (\tau) \right )
   \right ]
   \nonumber \\[1ex]
& = & e^{-\eta/\rho}\ \E \left [ \prod_{S \in \hyp_{1,\rho}}
   \left ( 1 - \prod_{\tau \in S \setminus \{ \rho \}} \theta_{m} (\tau) \right )
   \right ]\,,
   \label{eq:pois decomp}
\end{eqnarray}
where the first product is over hyperedges of cardinality up to~$k$
and the product over $\tau \in S \setminus \{ \rho \}$ is taken to be~1
if $S = \{ \rho \}$.
If $f : \Xi \to [0,1]$ is any function on a space $\Xi$ on which is
defined a Poisson process with intensity $\nu$, then the expected
product of $f$ at points of the Poisson process is given by
$$\exp \left [ \int (f(\xi) - 1) \, d\nu (\xi) \right ] \, .$$
Applying this to (\ref{eq:pois decomp}) with $\nu = \nu_\rho$
and $f (S) = 1 - \prod_{\tau \in S \setminus \{ \rho \}} \theta_{m} (\tau)$ gives
$$\log (1 - \theta_{m+1} (\rho)) = - \frac{\eta}{\rho}
   - \int \prod_{\tau \in S \setminus \{ \rho \}}
   \theta_{m} (\tau) \; d\nu (S) \, .$$
Break up the integral according to $|S|$.  Recall that for $j \geq 2$,
the law of $S \setminus \{ \rho \}$ on $\{ |S| = j \}$ is $\eta \mu_{j-1} / \rho$.
We may incorporate $-\eta / \rho$ as the $j=1$ term if we define
$\mu_0$ to be a point mass of~1 at the empty set and the
empty product to be~1.  These substitutions yield
$$\log (1 - \theta_{m+1} (\rho)) = - \frac{\eta}{\rho} \, \sum_{j'=0}^{k-1}
   \int \prod_{\tau \in S'} \theta_{m} (\tau) \; d\mu_{j'} (S) \, .$$
Here the primes are introduced to clarify the changes of
variable $j' = j-1, S' = S \setminus \{ \rho \}$.  We now drop the
primes and observe that $\mu_j$ is $1/(j!)$ times a product
measure.  Therefore the integral of the product factors, yielding
\begin{eqnarray*}
\log (1 - \theta_{m+1} (\rho)) & = & - \frac{\eta}{\rho} \, \sum_{j=0}^{k-1}
   \frac{1}{j!} \left ( \int_1^M \theta_{m} (\tau) \frac{d\tau}{\tau} \right )^j
   \\[1ex]
& = & - \frac{\eta}{\rho} \exp_k \left ( \int_1^M \theta_{m} (\tau) \frac{d\tau}{\tau}
   \right ) \, .
\end{eqnarray*}
Using the induction hypothesis again we substitute
$1 - e^{- \gamma_{m,M,k} (\eta) / \tau}$ for $\theta_{m} (\tau)$ to arrive at
$$\log (1 - \theta_{m+1} (\rho)) = - \frac{\eta}{\rho} \exp_k \left (
   \int_1^M \left ( 1 - e^{- \gamma_{m,M,k} (\eta) / \tau} \right ) \frac{d\tau}{\tau}
   \right ) \, .$$
Changing variables to $t = 1/\tau$ so that $dt / t = - d\tau / \tau$, yields
\begin{eqnarray*}
\log (1 - \theta_{m+1} (\rho)) & = & - \frac{\eta}{\rho} \exp \left (
   \int_{1/M}^1 \frac{1 - e^{- t \gamma_{m,M,k} (\eta)}}{t} \, dt
   \right ) \\
& = & - \frac{\eta}{\rho} A_M (\gamma_{m,M,k} (\eta)) \, .
\end{eqnarray*}
The right-hand side is equal to $- (1 / \rho) \gamma_{m+1,M,k} (\eta)$,
completing the induction.   $\Cox$

\begin{lemma} \label{lem:theta converge}
Fix any $\eta > \eta_*$.  Then
$$\theta_{m,k}^{M,\eta} (\rho) \to 1\,,$$
uniformly over $\rho$ in any bounded interval $[1,L]$
as $m , M , k \to \infty$.
\end{lemma}

\noindent{\bf Proof.}  The function
$z / \exp (A(z))$ is the real analytic function
$$
\exp \left( \int_1^z \frac{ e^{-u}}{u} \, du - \int_0^1 \frac{1 - e^{-u}}{u} \, du\right) = \exp (- \gamma - \Gamma (0,z)) ,
$$
where $\Gamma (0,z):=\int_z^\infty  e^{-t} \frac{dt}t$.  By (\ref{gamma}),
which evidently increases to $\eta_*$ as $z \uparrow \infty$.
It follows that for $\eta > \eta_*$, if we choose any positive
$\delta < (\eta / \eta_*) - 1$, then
$$\frac{\eta}{1 + \delta} > \eta_* > \frac{z}{e^{A(z)}}\,,$$
which implies that 
$$\eta e^{A(z)} > (1 + \delta) z\,,$$
for all $z > 0$.  Applying this to (\ref{eq:gamma}) with
$z = \gamma_{m,\infty,\infty}(u)$ leads to
$$\gamma_{m+1,\infty,\infty} (\eta) >
   (1 + \delta) \gamma_{m,\infty,\infty} (\eta)\,$$
which, in turn, leads inductively to
$$\gamma_{m,\infty,\infty} (\eta) > \eta_* (1 + \delta)^{m-1} \, .$$
Since $\gamma$ is increasing in all its arguments, this is true
for all greater $\eta$ as well.

Now, given $L , \ee > 0$, choose $m$ sufficiently large so that
$\gamma_{m,\infty,\infty} (\eta) > L \log(1/\ee)$.  The function
$\gamma$ is continuous in $M$ and $k$ at infinity, so we
may choose $M$ and $k$ such that $\gamma_{m,M,k} (\eta) > L \log (1/\ee)$.
It follows from Lemma~\ref{lem:theta} that
$$\theta_{m,k}^{M,\eta} (\rho) = 1 - e^{- \gamma_{m,M,k} (\eta) / \rho}
   > 1 - e^{- \log (1/\ee)} = 1 - \ee\,,$$
for $1 \leq \rho \leq L$, proving the lemma.   $\Cox$

\subsection{Proof of main theorems}
\label{ss:4.6}

\noindent{\bf Proof of Theorem~\ref{main_theorem}.}
Fix $\ee > 0$.  The first step is to use Lemma~\ref{lem:theta converge}
to pick $m,M,k$ such that
$$\theta_{m,k}^{M , \eta_* + \ee} (\rho) > \frac{3}{4}
   \;\;\; \mbox{ for all } \;\; 1 \leq \rho \leq L :=
   \exp \left ( \frac{3}{\ee} \right )
   \, . $$
Take $M$ to be larger if necessary so that we may assume $M \geq L$.
We deduce from the last displayed estimate with $\rho=p/y$ and from   Theorem~\ref{th:chi} that, for any prime $p$ in the interval $(y,My)$ and for $x$ sufficiently large, we have
$$\P_x \left ( \chi_{m,p}^{M , (\eta_* + \ee) J_0} \right ) >
   \frac{3}{4} \, .$$
Now let $Y$ be the number of $j$ in the interval
$I := [(\eta + \ee) J_0 , (\eta + 2 \ee) J_0]$ such that
$[X_j] = \{ p \}$ for some prime $p$ with $y < p < M y$
and $\chi_{m,p}^{M,j-1}$ holds.  Write $Y = \sum_{j \in I} Y_j$
where $Y_j$ is~1 if $[X_j] = \{ p \}$ for some prime $y < p < M y$
and zero otherwise.  We compute a lower bound on
$\E_x Y$ as follows.  The event $\chi_{m,p}^{M,j-1}$ is independent
of the event $[X_j] = \{ p \}$.  By (\ref{eq:unif})
and the definition of $J_0$ we have $\psi (x/p , y) / \psi(x,y) \sim (\log y) / p$.
Hence,
\begin{eqnarray*}
\E_x Y & = & \sum_{j \in I} \sum_{y < p < M y} \P ([X_j] = \{ p \} )
   \P (\chi_{m,p}^{M,j-1}) \\[1ex]
& = & \sum_{j \in I} \sum_{y < p < M y}
   \frac{\psi (x/p , y)}{x}
   \P (\chi_{m,p}^{M,j-1}) \\[1ex]
& \geq & \frac{1}{2} \sum_{j \in I} \sum_{y < p < M y}
   \frac{\pi (y)}{J_0} \frac{\log y}{p}
\end{eqnarray*}
for $x$ sufficiently large.  By the prime number theorem,
$$\sum_{y < p < M y} (\log y) / p \sim \log M \geq \log L
= 3 \ee^{-1} \, .$$
The outer sum has at least $\ee J_0$ terms, hence
\begin{equation} \label{eq:first moment}
\E_x Y \geq \frac{1}{2} (\ee J_0) \frac{\pi (y)}{J_0} (3 \ee^{-1})
   = \frac{3}{2} \pi (y) \, .
\end{equation}
In Lemma~\ref{lem:2mm} below, we will prove the second moment bound
$$\Cov (Y_i , Y_j) = o \left ( \frac{\pi (y)^2}{J_0^2} \right ) \, .$$
Using, this lemma,
\begin{eqnarray*}
\Var (Y) & = & \sum_{i , j \in I} \Cov (Y_i , Y_j) \\
& \leq & \E_x Y + 2 \sum_{i,j \in I, i < j} \Cov (Y_i , Y_j) \\[2ex]
& = & o(\pi (y)^2) \, .
\end{eqnarray*}
Together with (\ref{eq:first moment}), this implies that
$\P_x (Y > \pi (y)) \to 1$.  Recall from (\ref{eq:count})
that this implies more than $\pi (y)$ linear dependences
among the classes $[X_j]$ with $j \leq (\eta_* + 2 \ee) J_0$.
Since $\ee > 0$ was arbitrary, this completes the proof
of the theorem, modulo the lemma.   $\Cox$

\noindent{\bf Proof of Theorem~\ref{th:quant}.}
In the previous section, we chose $M$ to be absurdly large,
which allowed us to use only those $j$ in the interval
$[(\eta_* + \ee) J_0 , (\eta_* + 2 \ee) J_0]$.  We can
get much more reasonable values of $m,M$ and $k$ if we
are willing to let $\eta$ be a little bigger and to
use all the values of $j$ up to $\eta J$.  The computations
are in fact no harder (although the required convergence lemmas
did involve more work in the previous sections).

Fix $\eta , m , M$ and $k$ satisfying the inequality in the
hypothesis of the theorem.  Let
$$Z := \sum_{j=1}^J Z_J := \# \left \{ j \leq J : \chi_{m,k,p}^{M , j-1}
   \mbox{ occurs for all } p \in [X_j] \right \} \, .$$
Again, Lemma~\ref{lem:2mm} implies $\Var (Z) = o(\pi (y)^2)$.
If we are able to show
\begin{equation} \label{eq:1mm b}
\liminf_{x \to \infty} \frac{\E_x Z}{\pi (y)} > 1 \, ,
\end{equation}
then we would have $\P_x (Z > \pi (y)) \to 1$, which would imply
more than $\pi (y)$ linear dependences, thus establishing  the theorem.

To prove (\ref{eq:1mm b}), break
down $\E Z_j$ according to the value of $[X_j]$ and using
independence of $X_j$ from $\chi_{m,k,p}^{M,j-1}$.  This gives
\begin{eqnarray*}
\E_x Z_j & = & \sum_S \P_x ([X_j] = S)
   \P_x (\chi_{m,k,p}^{M,j-1}) \\[1ex]
& = & \sum_S \frac{\psi (x / \prod_{p \in S} p , y)}{x}
   \P_x (\chi_{m,k,p}^{M,j-1}) \\[1ex]
& \sim & \sum_S \frac{(\log y)^{|S|}}{\prod_{p \in S} p} \frac{\psi (x,y)}{x}
   \prod_{p \in S} \theta_{m,k}^{M,j/J_0} (p / y) \, .
\end{eqnarray*}
The final equality above used both equation (\ref{eq:unif})
and the formula (\ref{eq:chi square}) of Theorem~\ref{th:chi}.
Continuing, we use the identity $\psi (x,y) / x = \pi (y) / J_0$,
factor out this term, and rewrite the summand as a product:
$$\E Z_j = \frac{\pi (y)}{J_0} \sum_S \prod_{p \in S}
   \left ( \frac{\log y}{p} \theta_{m,k}^{M,j/J_0} (p/y) \right ) \, .$$
Let $B$ be any set and $\{ z_p : p \in B \}$ be any positive
real numbers with finite sum.  Let ${\cal B}$ denote the set
of finite subsets of $B$.  Then
$$\sum_{S \in {\cal B}} \prod_{p \in S} z_p
   = \prod_{p \in S} (1 + z_p) \to \exp (\sum_{p \in B} z_p)\,,$$
as $\max_{p \in B} z_p \to 0$.  Using this identity, we obtain
\begin{eqnarray*}
\E_x Z_j & \sim & \frac{\pi (y)}{J_0} \exp \left (
   \frac{1}{y} \sum_{y < p < M y} \frac{\log y}{p/y} \theta_{m,k}^{M,j/J_0}
   (p/y) \right ) \\[1ex]
& \sim & \frac{\pi (y)}{J_0} \exp \left ( \int_1^M \frac{1}{t}
   \theta_{m,k}^{M,j/J_0} (t) \, dt  \right )\,,
\end{eqnarray*}
by the prime number theorem.  The asymptotic equivalence is uniform
in $j \leq \eta J_0$.  Summing from $j=1$ to $\eta J_0$ now gives
\begin{eqnarray*}
\frac{\E_x Z}{\pi (y)} & \sim & \int_0^\eta \exp \left (
    \int_1^M \frac{1}{t} \theta_{m,k}^{M,u} (t) \, dt \right ) \; du \\[1ex]
& = & \int_0^\eta \frac{\gamma_{m+1} (u)}{u} \, du \, .
\end{eqnarray*}
By the hypothesized inequality, the right-hand side is greater than~1,
which establishes (\ref{eq:1mm b}) and completes the proof of the theorem.
$\Cox$
\begin{lemma} \label{lem:2mm}
Fix a finite real $M > 1$ and $\eta > 0$ and an integer $m \geq 1$.
Fix $1 \leq k \leq \infty$.  Then
$$\Cov (Z_i , Z_j) = o \left ( \frac{\pi (y)^2}{J_0^2} \right )$$
for all $1 \leq i < j \leq \eta J_0$.  The same is true with
$\Cov (Y_i , Y_j)$ in place of $\Cov (Z_i , Z_j)$.
\end{lemma}

\noindent{\bf Proof.}  Both arguments are the same, so we prove
this just for $\Cov (Z_i , Z_j)$.  It suffices to show that
$$\E_x (Z_i \cdot Z_j) \sim (\E_x Z_i) \cdot (\E_x Z_j) \, ,$$
uniformly for $1 \leq i < j \leq J$.  Conditioning on
$[X_i]$ and $[X_j]$, we see that this is the expectation of
$$\E_x(Z_i | [X_i] , [X_j]) \cdot \E_x (Z_j | [X_i] , [X_j]) \, . $$
The sets $[X_i]$ and $[X_j]$ are disjoint with probability going
to~1, so it suffices to show that
$\E_x(Z_i | [X_i] , [X_j])$ and $\E_x (Z_j | [X_i] , [X_j])$ are
asymptotically independent when $[X_i]$ and $[X_j]$ are
disjoint.  We have seen in Lemma~\ref{lem:conv 3} that the
collection of hypergraphs $\graph_{m,k,p}^{M,i-1,x}$ for
$p \in [X_i]$ and $\graph_{m,k,p}^{M,j-1,x}$ for $p \in [X_j]$
are disjoint and tree-like with probability going to~1,
and asymptotically independent.  The same is true of
the marked hypergraphs, by Lemma~\ref{lem:conv 4}.  Since
$Z_i$ is a bounded function of $[X_i]$ and the marked hypergraphs
$(\graph_{m,k,p}^{M,i-1,x} , \marked_{m,k,p}^{M,i-1,x})$ for
$p \in [X_i]$, and likewise for $[Z_j]$, we have the desired
conditional independence.
$\Cox$

\section{Implications for Factoring Algorithms}
\label{Algorithms}

In factoring algorithms we need to find 
a linear dependence mod 2 in our matrix of exponents.
We expect that the best algorithms known, due to Wiedemann or Lanczos
(see section 6.1.3 of  \cite{CP}), take  time
$$
\sim C \frac{y^2}{\log y\log\log y}
$$
for a positive constant $C$,
when we use the primes up to $y$ in our ``factor base''.
If we were to take $y=y_0$ then this number would be far larger than
$J_0$ and so would dominate the running time of the algorithm. Hence, to
optimize, we select $y=y_1$, which is far smaller, chosen to equalize
the running times of the two main parts of the algorithm, so that
\begin{equation} \label{optimize}
c \frac{\pi(y)}{\Psi(x,y)/x} \sim  \frac{y^2}{\log y\log\log y}
\end{equation}
for an appropriate constant $c>0$. One can show that one then has
$$
 y_1  =  y_0^{1 - (1+o(1))/\log\log x} ,
$$
with expected running time 
$$
  J_0\  y_0^{(1+o(1))/(\log\log x)^2}
$$
(see \cite{CGPT}).  

The proofs in the previous section work, as well for $y_1$, as for $y_0$.
In particular we can determine the speed-up for various choices of the 
parameters (though always with $m=\infty$, see \cite{CGPT} for more details):  
 
\vskip .25in
\centerline{\begin{tabular}{| l l | c c | c c | c c  | c c | }
\hline
&\ \ $k$ &&$M=\infty$  &&  $M=100$ &&  $M=10$ \\
\hline
&  $0$ &&  1 &&  1 && 1 \\
& $1$ &&  .7499 && .7517 && .7677\\
& $2$ &&  .6415 && .6448 && .6745\\
& $3$ &&  .5962 && .6011 && .6422\\
& $4$ &&  .5764  && .5823 && .6324\\
& $5$ &&  .567 && .575 && .630\\
\hline
\end{tabular}}
\bigskip
\centerline{\tt The value of $\eta$ such that there are $\sim \pi(y)$ }

\centerline{\tt pseudosmooths amongst the $a_j$ with   $j\leq \eta \pi(y)x/\Psi(x,y)$.}
\vskip .35in

So what effect will this reduction in the number of $a_j$ examined
have in the actual running time? Suppose that we replace $c$ in 
(\ref{optimize}) by $\eta c$, and  determine that the new running time
is given by (\ref{optimize}), after solving (\ref{optimize}) to determine $y=y_\eta$.

Now finding this solution is tantamount to finding a 
solution to $h(u_\eta)=\log (c\eta\log\log y)$ where
$h(u):= \frac 1u \log x+\log \rho(u)$. We have
$h'(u)= -1-(1+o(1))/\log u)$ and so 
$u_1-u_\eta=\log \eta  (1-(1+o(1))/\log u)$.
Our running time  therefore changes by a factor of
\begin{eqnarray*}
\sim   x^{\frac 2{u_\eta}-\frac 2{u_1}} & = & 
\exp \left( \frac {2(u_1-u_\eta)\log x}{u_1u_\eta} \right) =
\exp \left( \frac {2\log \eta \log x}{u_1^2} \left(1-\frac{1+o(1)}{\log u} \right) \right)
\\ & = & \exp \left(   \log \eta ( \log\log x   + \log\log\log x-\log 2-4+o(1)) \right) \\
& = &  \left( \frac {2e^4+o(1)}{ \log x \log\log x }   \right)^{ \log (1/\eta)}\,,
\end{eqnarray*}
since $\log^2 y_1 = \log^2 L(x) \left( 1 + \frac{\log\log\log x-\log 2-4+o(1)}{\log\log x}  \right)$.

\

Data on the effect of large prime variations that has been gathered from running factoring algorithms, 
seems rather  different from what we have
obtained here. One  reason  for this is that, in our analysis,
the variations in $M$ and $k$ simply affect the number of $a_j$ being considered,
whereas in reality these affect
not only the number of $a_j$ being considered, but also several other important
quantities. For instance,  the amount of sieving that needs to be done, and
also the amount of
data that needs to be ``swapped'' (typically one saves the $a_j$ with several
large prime factors to the disk, or somewhere else suitable for a lot of data).
It is an interesting problem to try to properly analyze the construction of
programs, so as to incorporate the results that we have obtained and to get
predictions that would help the choice of parameters in computer algorithms.

\end{document}